\documentclass[notitlepage,11pt,reqno]{amsart}
\usepackage[foot]{amsaddr}
\usepackage[normalem]{ulem}
\usepackage{array,amssymb,nicefrac,bm,upgreek,mathtools,verbatim,comment}
\usepackage[mathscr]{eucal}
\usepackage{dsfont}
\usepackage[final]{hyperref}
\usepackage{amsopn}

\usepackage[margin=1in]{geometry}
\allowdisplaybreaks
\hypersetup{
colorlinks=true,
citecolor=mblue,
linkcolor=mblue,
urlcolor = blue,
anchorcolor = blue,
breaklinks}
\usepackage[msc-links,abbrev]{amsrefs}
\usepackage{color}
\definecolor{dmagenta}{rgb}{.4,.1,.5}
\definecolor{dblue}{rgb}{.0,.0,.5}
\definecolor{mblue}{rgb}{.0,.0,.7}
\definecolor{ddblue}{rgb}{.0,.0,.4}
\definecolor{dred}{rgb}{.7,.0,.0}
\definecolor{dgreen}{rgb}{.0,.5,.0}
\definecolor{Eeom}{rgb}{.0,.0,.5}
\definecolor{cm}{cmyk}{1,.0,.0,.0}

\numberwithin{equation}{section}

\theoremstyle{plain}
\newtheorem{theorem}{Theorem}[section]

\theoremstyle{definition}
\newtheorem{assumption}{Assumption}[section]
\newtheorem{definition}{Definition}[section]

\newtheorem{algorithm}{Algorithm}[section]

\theoremstyle{remark}

\usepackage[capitalize,nameinlink]{cleveref}

\crefname{section}{Section}{Sections}
\crefname{subsection}{Subsection}{Subsections}
\crefname{condition}{Condition}{Conditions}
\crefname{hypothesis}{Hypothesis}{Conditions}
\crefname{assumption}{Assumption}{Assumptions}
\crefname{lemma}{Lemma}{Lemmas}
\crefname{definition}{Definition}{Definitions}
\crefname{fact}{Fact}{Facts}
\crefname{corollary}{Corollary}{Corollaries}
\crefname{theorem}{Theorem}{Theorems}
\crefname{claim}{Claim}{Claims}

\Crefname{figure}{Figure}{Figures}

\crefformat{equation}{\textup{#2(#1)#3}}
\crefrangeformat{equation}{\textup{#3(#1)#4--#5(#2)#6}}
\crefmultiformat{equation}{\textup{#2(#1)#3}}{ and \textup{#2(#1)#3}}
{, \textup{#2(#1)#3}}{, and \textup{#2(#1)#3}}
\crefrangemultiformat{equation}{\textup{#3(#1)#4--#5(#2)#6}}%
{ and \textup{#3(#1)#4--#5(#2)#6}}{, \textup{#3(#1)#4--#5(#2)#6}}%
{, and \textup{#3(#1)#4--#5(#2)#6}}

\Crefformat{equation}{#2Equation~\textup{(#1)}#3}
\Crefrangeformat{equation}{Equations~\textup{#3(#1)#4--#5(#2)#6}}
\Crefmultiformat{equation}{Equations~\textup{#2(#1)#3}}{ and \textup{#2(#1)#3}}
{, \textup{#2(#1)#3}}{, and \textup{#2(#1)#3}}
\Crefrangemultiformat{equation}{Equations~\textup{#3(#1)#4--#5(#2)#6}}%
{ and \textup{#3(#1)#4--#5(#2)#6}}{, \textup{#3(#1)#4--#5(#2)#6}}%
{, and \textup{#3(#1)#4--#5(#2)#6}}

\crefdefaultlabelformat{#2\textup{#1}#3}

\newcommand{\Ind}{\mathds{1}}       

\newcommand{\Act}{{\mathds{U}}}
\newcommand{\cA}{\mathcal{A}}
\newcommand{\uuptau}{{\Breve{\uptau}}}
\newcommand{\sB}{{\mathscr{B}}}
\newcommand{\cB}{{\mathcal{B}}}     
\newcommand{\Cc}{{C}}               
\newcommand{\sD}{{\mathscr{D}}}
\newcommand{\sE}{{\mathscr{E}}}     
\newcommand{\sF}{\mathfrak{F}}      
\newcommand{\cG}{{\mathcal{G}}}     
\newcommand{\sG}{{\mathscr{G}}}     
\newcommand{\sH}{{\mathscr{H}}}
\newcommand{\cI}{{\mathcal{I}}}
\newcommand{\sK}{{\mathscr{K}}}     
\newcommand{\cK}{{\mathcal{K}}}
\newcommand{\sL}{{\mathscr{L}}}
\newcommand{\Pm}{{\mathcal{P}}}     
\newcommand{\cS}{{\mathcal{S}}}     
\newcommand{\cT}{{\mathcal{T}}}

\newcommand{\Lyap}{{\mathscr{V}}}   
\newcommand{\cX}{{\mathcal{X}}}

\newcommand{\Uadm}{\mathfrak{U}}
\newcommand{\Um}{\mathfrak{U}_{\mathrm{m}}}
\newcommand{\Usm}{\mathfrak{U}_{\mathrm{sm}}}

\newcommand{\Uc}{\mathfrak{U}_{\mathrm{causal}}}
\newcommand{\lamstrdm}{\lambda^{*,{\rm m}}}
\newcommand{\lamstrdf}{\lambda^{*,{\rm d}}}
\newcommand{\lamstrcm}{\lambda^{*,{\rm c}}}
\newcommand{\PW}{\mathcal{M}(\widetilde{W})}

\newcommand{\RR}{\mathds{R}} 
\newcommand{\NN}{\mathds{N}} 
\newcommand{\Rd}{{\mathds{R}^{d}}}
\DeclareMathOperator{\Exp}{\mathbb{E}} 
\DeclareMathOperator{\Prob}{\mathbb{P}} 
\newcommand{\D}{\mathrm{d}} 

\newcommand{\abs}[1]{\lvert#1\rvert}
\newcommand{\norm}[1]{\lVert#1\rVert}

\newcommand{\transp}{^{\mathsf{T}}}

\DeclareMathOperator*{\Argmin}{Arg\,min}
\DeclareMathOperator*{\Argmax}{Arg\,max}

\DeclareMathOperator*{\diag}{diag}
\DeclareMathOperator*{\trace}{trace}

\DeclareMathOperator*{\sgn}{sign}

\newcommand{\grad}{\nabla}

\newcommand{\df}{\coloneqq}
\newcommand{\sorder}{{\mathfrak{o}}} 

\makeatletter
\DeclareRobustCommand\widecheck[1]{{\mathpalette\@widecheck{#1}}}
\def\@widecheck#1#2{%
    \setbox\z@\hbox{\m@th$#1#2$}%
    \setbox\tw@\hbox{\m@th$#1%
       \widehat{%
          \vrule\@width\z@\@height\ht\z@
          \vrule\@height\z@\@width\wd\z@}$}%
    \dp\tw@-\ht\z@
    \@tempdima\ht\z@ \advance\@tempdima2\ht\tw@ \divide\@tempdima\thr@@
    \setbox\tw@\hbox{%
       \raise\@tempdima\hbox{\scalebox{1}[-1]{\lower\@tempdima\box
\tw@}}}%
    {\ooalign{\box\tw@ \cr \box\z@}}}
\makeatother

\usepackage{accents}
\newlength{\dhatheight}

\begin{document}
\title[A survey on ergodic risk-sensitive control]%
 {Ergodic Risk-sensitive control---A survey}

\author[Anup Biswas]{Anup Biswas$^\dag$}
\address{$^\dag$ Department of Mathematics\\
Indian Institute of Science Education and Research Pune\\
Dr.\ Homi Bhabha Road\\
Pune 411008, India}
\email{anup@iiserpune.ac.in}

\author[Vivek S. Borkar]{Vivek S. Borkar$^\ddag$}
\address{Department of Electrical Engineering,
Indian Institute of Technology,
Powai, Mumbai 400076, India}
\email{borkar@ee.iitb.ac.in}


\date{\today}

\begin{abstract}
Risk-sensitive control has received considerable interest since the seminal 
work of Howard and Matheson \cite{Howard-71}  because of its ability to account for fluctuations about the mean, its connection with $H_\infty$ control, and its application to financial mathematics. In this article we attempt to put together
a comprehensive survey on the research  done on ergodic risk-sensitive 
control over the last four decades.
\end{abstract}


\maketitle
\tableofcontents

\section{Introduction}

Given a controlled stochastic process $\textbf{X}=\{X_t\}$ on a
state space $\cS$, controlled by 
the process $\zeta$, the ergodic risk-sensitive cost is defined as
$$\sE_x(c, \zeta)\df\limsup_{T\to\infty} \frac{1}{\gamma T}
\log \Exp_x\left[e^{\int_0^T \gamma c(X_t, \zeta_t)\D{t}}\right],\quad x\in\cS,$$
$c$ being the running cost and $\gamma\neq 0$ being the risk-parameter. The 
ergodic risk-sensitive control (ERSC) problem is about studying the
minimization problem
$$\lambda^*=\inf_{x\in\cS}\, \inf_{\zeta}\sE_x(c, \zeta).$$
Replacing `inf' with `sup' leads to the corresponding reward maximization problem which we discuss briefly later. Observe , however, that unlike the classical cost functionals, the reward maximization problem is \textit{not} equivalent to the cost minimization problem obtained by flipping the sign of the instantaneous reward. 

 Suppose $c$ is non-negative (more generally, bounded from below). The decision maker is supposed to be risk-averse or risk-sensitive for
$\gamma>0$, risk-neutral for $\gamma=0$ and risk-seeking for $\gamma<0$. 
The risk-neutral case, in a suitable limiting sense,  corresponds to the the classical ergodic control problem which has already been studied extensively
(see \cites{ABFGM93,red-book} and references therein). 
The goal of this article is to review the development 
of ERSC problems when $\gamma\neq 0$. The study of ERSC can be traced
back to the seminal work of Howard and Matheson \cite{Howard-71} where 
the problem was studied for controlled Markov chains with finite state and action sets. Since then this
area has been developed intensively in the past forty years. A major fillip came from a series of works by Peter Whittle in the eighties, culminating in \cite{Whit-90}. One major motivation was the strongly felt need for criteria going beyond those based purely on mean rewards, that did not put any weight whatsoever on fluctuations around the mean. The obvious extensions such as considering a weighted sum of mean and variance in some form (the Markowitz model in finance being the prime example) faced problems such as non-availability of the `principle of time-consistency' or the dynamic programming principle. As the exponential function can be viewed as the weighted sum of all powers, its expectation is a weighted sum of all moments. Thus its expctation does account for higher moments. In addition, by facilitating a multiplicative form of dynamic programming (as opposed  to the additive form for classical criteria), it does obey the principle of time consistency. This made risk-sensitive control an attractive proposition.

Two classical applications of ergodic risk-sensitive control problem motivated by such considerations  come
from {\it robust control theory} and {\it portfolio optimization problems}.
\begin{itemize}
\item[--] ({\it Robust control theory}) Since it is often almost impossible to find a {\it true model} of a system,  robust control theory seeks criteria that could deal with model uncertainty.
The connection between risk-sensitive control and robust control started with the work of Glover and Doyle \cite{MR960663} (see also Whittle \cites{Whit-81,Whit-90}). Risk-sensitive minimization problems naturally give rise to two person zero-sum
differential games (first found in the work of Jacobson \cite{J-73}) which are of interest in the robust control theory. In the differential game formulation there are two players, one representing the disturbance entering the system which will attempt to degrade system performance, and the other representing the {\it actual} control for the system. Readers may also consult \cite{Dupuis2000} for more on {\it power gain inequality} and its connection to the ergodic risk-sensitive value. We shall briefly discuss this and its connection to $H_\infty$ control  in Section~\ref{H-infinity}.

Another way to deal with the model uncertainty is to consider partially observed or Hidden Markov chain models. ERSC has been studied in these
frameworks as well. It should be note that under fairly general assumptions
and suitable change of measures, partially observed models can be changed into a fully observed control problems (cf. \cite{MR1276773}). Interested readers
may consult \cite{MR1276773,Fleming-99,MR1427787,MR1690565,BJ97} for more details in this direction.
\\[2mm]

\item[--] ({\it Portfolio optimization}) The risk-sensitive formulation of the portfolio optimization problem was introduced in the seminal works of 
Bielecki and Pliska \cite{MR1675114}, Fleming and Sheu \cite{MR1722286}. Since then this area has grown substantially (see \cites{MR1802598,MR1910647,MR1890061,MR1972534,MR2193508} and references therein). In addition to the considerations already discussed, multiplicative / exponential models arise naturally in finance due to `compounding' effects. Suppose that there
are $N$ risky assets and the investor allocates fraction $u^i_t$ of its wealth to the $i$-th risky asset, $i=1, 2, \ldots, N$. The total wealth $V_t$, at time $t$, of the investor 
is then given by
\begin{equation}\label{E1.1}
\D{V}_t= V_t\left[ r_t \left(1-\sum_{i=1}^N u^i_t\right) \D{t}+\sum_{i=1}^N u^i_t \frac{\D{S}^i_t}{S^i_t}\right],
\end{equation}
where $r$ denotes the risk-free interest rate, $S^i$ the share price of the $i$-th risky asset. Let $\Act\subset\RR^N$ be a constrain set and $\zeta_t=(u^1_t, \ldots, u^N_t)\in \Act$ for all $t$.
In portfolio optimization, one wishes to maximize the long term value of $\gamma^{-1}\Exp[V^\gamma_T]$, for some $\gamma\in (-\infty, 1)\setminus\{0\}$, over all possible
investment allocations. Now suppose that there are $d$ economic factors given by the vector $\tilde{X}_t\in\Rd$ that governs the market performance and evolves according to
the stochastic differential equation
$$\D{\tilde{X}_t}=\tilde{b}(\tilde{X}_t)\D{t} + \D{W}_t,$$
where $W$ is a $d$-dimensional standard Brownian motion. The share price dynamics is given by
$$\frac{\D{S}^i_t}{S^i_t}= \mu^i(\tilde{X}_t) \D{t} + \upsigma^i_D \cdot \D{W}_t + \upsigma^i_I \cdot \D{\tilde{W}}_t, \quad i=1, \ldots, N, $$
where $\tilde{W}$ is an $N$ -dimensional standard Brownian motion independent of $W$. Assume that $\upsigma^i_D, \upsigma^i_I, r$ are constant vectors. Applying It\^{o}'s formula one 
can easily find from \eqref{E1.1} the differential equation satisfied by $\log V_t$.
Then defining
\begin{gather*}
b(x, u)\df \tilde{b}(x) + \gamma\sum_{i=1}^N u^i \upsigma^i_D, \quad \bar\mu_i(x)\df \mu^i(x)-r, \quad 
\upsigma^i=[\upsigma^i_D, \upsigma^i_I]\in\RR^{d+N},
\\
\ell(x, u)\df-\frac{1}{2}(1-\gamma)\abs{\sum_{i=1}^N u^i \upsigma^i}^2 + \sum_{i=1}^N u^i \bar\mu^i(x) + r,
\end{gather*}
one can check that above maximization problem is equivalent to maximizing (see \cite{MR1910647})
$$\log\Exp\left[e^{\gamma \int_0^T \ell(X_t, \zeta_t)\D{t}}\right],$$
where $\zeta_t\in\Act$ and 
$$\D{X}_t=b(X_t, \zeta_t)\D{t} + \D{W}_t.$$
Thus the long-term asymptotics (that is, as $T\to\infty$) corresponds to the ergodic risk-sensitive control problems.

It is worth noting that some of the early work in this direction came from information theorists, notably Thomas Cover and his associates. See \cite{Cover}, Chapter 16, and its bibliographical note.
\end{itemize}

We mention in passing another application, viz.\ to minimizing or maximizing the asymptotic rate of exit of a controlled Markov process from a prescribed subset of its state space. This can be reduced to a risk-sensitive control problem \cite{BBEigen}, \cite{BFilar}.

The main focus of this article is on ERSC problems. There is also an enormous amount of 
work done on the finite horizon version of risk-sensitive control problems which we do not discuss in this article. Interested readers
may look at \cites{Bensoussan-95a,MR3286649,MR3625748,MR3173005,Bensoussan-03,Bensoussan=-95b,Bensoussan-96,Jaskiewicz-14,Whit-81}.
Broadly speaking, the ERSC problems are treated in three different ways.
The first one corresponds to the variational representation of the moment
generating function. This helps us to transform the above minimization problem to an ergodic zero-sum game problem (see \cites{MR3900793,Fleming-95a,Fleming-91,Fleming-97a}). 
 The second approach for solving 
ERSC problem is an approximation method based on the discounted risk-sensitive problem
(see \cites{DiMasi-99b,MR2174017,MR2799396}). Discounted risk-sensitive control is not amenable to dynamic programming, but by treating the risk-snsitivity parameter as a variable, one manages to make the problem analytically tractable. The dynamic programming equation of the risk-sensitive control problem is a 
nonlinear eigenvalue problem. The third approach is
more direct where the nonlinear eigenvalue problem is analyzed using
Krein-Rutman theorem (see \cites{MR3926044,MR2808061,MR3629428}). We divide the review of ERSC problems
 in three major parts, namely, discrete time set up, controlled diffusions and  continuous time Markov chains, wherein we touch upon all three approaches above.

\vspace{.2in}

The following is a list of the abbreviations used in this paper
\[
\begin{array}{r l}
\text{DTCMC} & \text{discrete time controlled Markov chain} \\
\text{ERSC} & \text{ergodic risk-sensitive control}\\
\text{CTCMC} & \text{continuous time controlled Markov chain}\\
\text{PIA} & \text{policy iteration algorithm}\\
\text{RVI} & \text{relative value iteration}
\end{array}
\hspace{1000pt minus 1fill} 
\]
We also summarize key notations used in this article
\[
\begin{array}{r l}
\sB(\cX) & \text{Borel $\sigma$ algebra on the topological space $\cX$}\\
C_b(\cX) & \text{set of all real-valued bounded, continuous functions on $\cX$}\\
C^k(\cX) & \text{set of all $k$-times continuous differentiable functions on $\cX\subset\Rd$}\\
C^k_{+}(\cX) & \text{subset of functions of $C^k(\cX)$ that are positive on $\cX$}\\
\lambda^{*,{\rm m}} & \text{optimal ergodic risk-sensitive value for DTCMC}\\
\lambda^{*,{\rm d}} & \text{optimal ergodic risk-sensitive value for controlled diffusion}\\
\lambda^{*,{\rm c}} & \text{optimal ergodic risk-sensitive value for CTCMC}
\end{array}
\hspace{1000pt minus 1fill} 
\]

\section{Risk-sensitive control of Discrete time Markov chains}\label{S-DTMC}
We begin by introducing the general setting of a controlled discrete time Markov chain.
Consider a controlled Markov process $\textbf{X} \df \{X_0, X_1,\dots\}$ on a Borel space
$\cS$
controlled by a control process $\zeta \df \{\zeta_0,\zeta_1,\dots \}$ taking values in $\Act$. Here $\Act$ is a Borel space endowed with the Borel $\sigma$ algebra $\sB(\Act)$. For every $x\in \cS$, $\Act(x) \in \sB(\Act)$ stands for the nonempty compact set of all admissible actions when the system is at the state $x$. The space of all admissible state action pairs is given by $\sK \df \{(x,u) : x\in S, u\in \Act(x)\}$. For each $A\in \sB(\cS)$ the controlled stochastic kernel $P(A |\cdot):\sK\to [0,1]$ is Borel measurable. 
We denote by $c:\sK \to \RR_{+}$ the one-stage cost function. For each $t\in\NN$, the space $\sH_t$ denotes the {\it admissible histories} up to time $t$, where $\sH_0 := \cS$, $\sH_t = \sK\times\sH_{t-1}$.
A generic element $h_t$ of $\sH_t$ is a vector of the form
$$
h_t=(x_0, u_0, x_1, u_1, \ldots, x_{t-1}, u_{t-1}, x_t),
\quad \text{with}\;\;(x_s, u_s)\in\sK,\quad 0\leq s\leq t-1, \ x_0\in \cS,
$$
denotes the observable history of the process up to time $t$.
Let us also denote by $\sF_n=\sB(\sH_n) :=$ the Borel $\sigma$-field of $\sH_n$.
An {\it admissible control}  is a sequence $\zeta = \{\zeta_0,\zeta_1,\dots \}$ where for each $t\in\NN$\,, $\zeta_t : \sH_t\to\Act$ is a measurable map satisfying $\zeta_t(h_t)\in\Act({x_t})$, for all $h_t\in\sH_t$. The set of all admissible policies is denoted by $\Uadm$. It is well known that for a given initial state $x\in S$ and policy $\zeta\in\Uadm$ there exists a unique probability measure $\Prob_x^{\zeta}$ on 
$(\Omega, \sB(\Omega))$, where $\Omega=(S\times\Act)^\infty$, (see 
\cite[p.4]{H89}, \cite{ABFGM93}) satisfying the following
\begin{equation}\label{Markov1}
\Prob_x^\zeta(X_0=x)=1,\quad \text{and}\quad
\Prob_x^{\zeta}(X_{t+1} \in A | \sH_t, \zeta_t) = P(A | X_t, \zeta_t)\quad \forall \,\, A\in \sB(S)\,.
\end{equation}
The corresponding expectation operator is denoted by $\Exp_x^{\zeta}$. A policy $\zeta\in\Uadm$ is said to be a Markov policy if $\zeta_t(h_t) = v_t(x_t)$ for all $h_t\in\sH_t$, for some measurable map 
$v_t:\cS\to\Act$ such that $v_t(x)\in\Act(x)$ for all $x\in\cS$. The set of all Markov policies is denoted by $\Um$. If the map $v_t$ does not have any explicit time dependence, that is, $\zeta_t (h_t) = v(x_t)$ for all $h_t\in\sH_t$, then $\zeta$ is called a stationary Markov strategy and  we denote the 
set of all stationary Markov strategies by $\Usm$.(We use the words `strategy' and `policy' interchangeably.)  From \cite[p.6]{H89} (also see \cite{ABFGM93}), it is easy to see that under any Markov policy $\zeta\in\Um$, the corresponding stochastic process $\textbf{X}$ is strong Markov. For each $\zeta\in\Uadm$, the ergodic risk-sensitive cost is given by
\begin{equation}\label{EErgocost}
\sE_x(c, \zeta) \,\df\, \limsup_{ T\to\infty} \, \frac{1}{\gamma T}\,
\log \Exp_x^{\zeta} \left[e^{\sum_{t = 0}^{T-1} \gamma c(X_t, \zeta_t)}\right],
\end{equation}
where $\gamma\neq 0$ and $\textbf{X}$ is the discrete time controlled Markov chain (DTCMC)
corresponding to the control $\zeta\in\Uadm$, with initial state $x$.
Our aim is to minimize \cref{EErgocost} over all admissible policies $\Uadm$. In other words,
we are interested in the quantity
\begin{equation}\label{lamstr-DMC}
\lamstrdm= \, \inf_{x\in \cS}\, \inf_{\zeta\in \Uadm}\sE_x(c, \zeta).
\end{equation}
We refer to this as an ergodic risk-sensitive control (ERSC) problem.
A policy $\zeta^{*}\in \Uadm$ is said to be optimal if for all $x\in\cS$
$$\sE_x(c,\zeta^{*}) \, = \, \inf_{x\in S}\,\inf_{\zeta\in \Uadm}\sE_x(c, \zeta).$$
Note that in general,  $\sE_x(c,\zeta)$ is not independent of 
$x$ for $\zeta\in\Usm$. Let us also mention the optimality 
equation which will be important for characterizing the optimal stationary
Markov controls.
\begin{definition}\label{A2.1}
A positive function $\psi:\cS\to (0, \infty)$ and a real number
$\lambda$ are said to form an eigen-pair $(\psi, \lambda)$ if
\begin{equation}\label{ED2.1}
\sgn(\gamma)e^{\gamma\lambda}\psi(x) = \min_{u\in\Act(x)}\left[\sgn(\gamma)e^{\gamma c(x,u)}\int_{\cS} \psi(y)P(\D{y}|x,u)\right]\quad\text{for}\,\, x\in \cS.
\end{equation}
We call $\psi$ an eigenfunction corresponding to the eigenvalue
$\gamma\lambda$.
\end{definition}
We impose the following standard assumption on our model.
\begin{assumption}
The following hold.
\begin{itemize}
\item[(i)] The transition kernel $P(\cdot| x, u)$ is weakly continuous in
$(x, a)$, that is, for every $f\in\Cc_b(\cS)$ we have 
$\int_{\cS} f(y)P(\D{y}|x, u)$ continuous in $\sK$.

\item[(ii)] $u\mapsto c(x, u)$ is continuous in $\Act(x)$ for all $x\in\cS$.
\end{itemize}
\end{assumption}

\subsection{Finite state space}
Suppose that $\cS$ is a finite set. The very first  ERSC control problem appeared in the work of Howard and Matheson \cite{Howard-71} where the authors studied an ergodic risk-reward problem
under the assumption that $\textbf{X}$ is irreducible and aperiodic under every stationary Markov policy. 
Since then the finite state situation has been studied in several works 
\cite{BCF98,MR3639105,MR2115041,MR2510639,MR2799396,MR2558431,MR2115041,MR1893292,Fleming-99,MR736636,MR2932929,MR1690565}. For instance, for a
(uncontrolled) Markov chain $\textbf{X}$ with transition matrix $P$, it is well-known that
$$x\mapsto \limsup_{ T\to\infty} \, \frac{1}{\gamma T}\,
\log \Exp_x \left[e^{\sum_{t = 0}^{T-1} \gamma c(X_t)}\right]$$
is constant on each communicating class (cf. \cite[Lemma~1]{BCF98}). Moreover, if  $\textbf{X}$ is irreducible, then 
$$\limsup_{ T\to\infty} \, \frac{1}{\gamma T}\,
\log \Exp_x \left[e^{\sum_{t = 0}^{T-1} \gamma c(X_t)}\right]= \frac{1}{\gamma} \log \rho(\tilde{P})$$
where $\tilde{P}_{ij}\df P_{ij} e^{\gamma c(i)}$ and $ \rho(\tilde{P})$ denotes the spectral radius of $\tilde{P}$. Furthermore, the existence of such an eigen-pair can be characterized by the following result.
\begin{theorem}[\cite{MR2510639}]
Let $P_{xy}=P(y|x)$ for $x, y\in\cS$. Then the following are equivalent.
\begin{itemize}
\item[(i)] For each $c:\cS\to \RR$ there exists an eigen-pair $(\psi, \lambda)$ satisfying
$$e^{\gamma\lambda} \psi(x)= e^{\gamma c(x)} \Exp_x[\psi(X_1)]=e^{\gamma c(x)} \sum_{y\in\cS} \psi(y) P_{xy}, \quad \text{for all}\; x\in\cS.$$ 
\item[(ii)] For every cost function $c$ the mapping
$$x\mapsto \limsup_{ T\to\infty} \, \frac{1}{\gamma T}\,
\log \Exp_x \left[e^{\sum_{t = 0}^{T-1} \gamma c(X_t)}\right]$$
is constant.
\item[(iii)] The transition matrix $P$ has a unique recurrent class $\mathcal{C}\subset\cS$ and there exists a constant $m$ such that
$$\Prob_x(\uuptau_{\mathcal{C}}\leq m)=1\quad \text{for all}\; x\in\cS$$
\end{itemize}
where $\uuptau_{\mathcal{C}}$ denotes the return time to the set $\mathcal{C}$, that is, 
$$\uuptau_{\mathcal{C}}\df \inf\{n\geq 1\; :\; X_n\in\mathcal{C}\}.$$
\end{theorem}
Moreover, if one of the above conditions holds, then the eigenfunction $\psi$ can be represented as follows \cite{MR2510639,MR1893292}
\begin{equation}\label{E2.5}
\psi(x)=\Exp_x\left[e^{\gamma \sum_{i=0}^{\uuptau_z-1} (c(X_i)-\lambda)}\right] \quad \forall \ x\in \cS\setminus\{z\}, \ \psi(z)=1,
\end{equation}
where $\uuptau_z=\uuptau_{\{z\}}$. $\lambda$ is the value of the average risk-sensitive cost. This representation of the eigenfunction will be crucial in our study and will appear in several places below. The following result on the ERSC problems can be found in \cite[Theorem~3.1]{MR1893292} (see also \cite{MR1687362,MR1466928})
\begin{theorem}\label{T2.2}
Suppose that under every stationary policy $\zeta\in\Usm$, each pair of states in $\cS$ communicates under $\textbf{X}$. Then the following hold for every $\gamma\neq 0$.
\begin{itemize}
\item[(i)] There exists an eigen-pair $(\Psi, \lambda), \Psi>0,$ satisfying
\begin{equation}\label{ET2.2A}
\sgn(\gamma)e^{\gamma\lambda}\Psi(x) = \min_{u\in\Act(x)}\left[\sgn(\gamma)e^{\gamma c(x,u)}\sum_{y\in\cS} \Psi(x)P(y|x,u)\right]\quad\text{for}\,\, x\in \cS.
\end{equation}
\item[(ii)] $\inf_{\zeta\in \Uadm}\sE_x(c, \zeta)=\lamstrdm=\lambda$ for each $x\in\cS$.
\item[(iii)] Every minimizing selector of \eqref{ET2.2A} is an optimal policy.
\item[(iv)] $(\Psi, \lambda)$ satisfying \eqref{ET2.2A} is unique provided we set $\Psi(z)=1$ for a prescribed state $z\in\cS$.
\end{itemize}
\end{theorem}
Note that the above result requires the DTCMC to be communicating under every stationary policy. Theorem~\ref{T2.2} also appears in \cite{MR1732397}
where it is proved under an additional assumption that $P(x|x, u)>0$ for all $(x, u)\in\sK$.
 In \cite{MR1973378} the author shows that given any two states $x, y\in\cS$, if we can find a 
stationary policy under which $y$ is accessible from $x$, then there exists
$\Lambda_0>0$ such that an eigen-pair satisfying \eqref{ET2.2A}
exists for $\gamma$ satisfying
$\gamma \norm{c}_{\rm sp}<\Lambda_0$ ($\norm{\cdot}_{\rm sp}$ denotes the span  semi-norm  defined as $\norm{c}_{\rm sp}=\sup_{x, u} c- \inf_{x, u} c$). Also, note that the hypothesis of a single communicating class for
every stationary control is important to ensure that 
$\inf_{\zeta\in \Uadm}\sE_x(c, \zeta)$ is independent of $x$
(a specific example can be found in \cite[Proposition~3.1]{MR1687362}).

Consider the assumption:
\begin{assumption}[Simultaneous Doeblin Condition]\label{A2.2}
There exists a state $z\in\cS$ and a positive integer $K$ such that
$$\Exp_x^\zeta[\uuptau_z]\leq K\quad \text{for all}\; x\in\cS, \quad \text{and}\; \zeta\in\Usm.$$
\end{assumption}
A general characterization of the optimal value is then obtained in
\cite[Theorem~3.5]{MR2115041}.

\begin{theorem}
Suppose that $\Act$ is a finite set and Assumption~\ref{A2.2} holds.
Then for every $x\in\cS$ we have
$$\inf_{\zeta\in \Usm}\sE_x(c, \zeta)=\inf_{g\in\mathfrak{G}}g(x),$$
where $\mathfrak{G}$ denotes the collection all functions
$g:\cS\to\RR$ satisfying
\begin{itemize}
\item[(i)] For each $x\in\cS$
$$g(x)=\min_{u\in\Act(x)}\, \max\{g(y)\; :\; P(y|x, u)>0\}.$$
\item[(ii)] There exists a positive function $h$ such that
$$e^{\gamma g(x)}h(x)\geq \min_{\mathcal{B}_g(x)}\left[
e^{\gamma c(x, u)}\sum_{y\in\cS} h(y) P(y|x, u)\right]
\quad x\in\cS,$$
where
$$\mathcal{B}_g(x)\df
\Bigl\{u\in\Act(x)\; :\; g(x)=\max\{g(y)\, :\, P(y|x, u)>0\}\Bigr\}.$$
\end{itemize}
\end{theorem}
A generalization of the above result for DTCMC with a general state space can be found in \cite{MR2585141}.

\subsection{Countable state space}
Now suppose that $\cS$ is countable. Without any loss of generality,
assume that $\cS=\{0, 1, 2,\ldots\}$.
The analysis of ERSC problem becomes
more involved due to non-compactness of $\cS$. If the running
cost $c$ is bounded, then a result analogous to Theorem~\ref{T2.2}
is possible, provided $\gamma$ is small.

\begin{theorem}[\cite{MR1687362}]
Let Assumption~\ref{A2.2} hold. Define
$$
\mu=\frac{\log(K+1)-\log K}{K+1}, 
\quad \norm{c}=\sup_{\sK}|c(x, u)|.
$$
Then for each $0\neq \gamma\in (-\frac{\mu}{2\norm{c}}, \frac{\mu}{2\norm{c}})$ there exists an eigen-pair $(\Psi, \lambda)$
with bounded  $\Psi$ that satisfies \eqref{ET2.2A}. Furthermore,
the conclusions of Theorem~\ref{T2.2} (ii)-(iv) hold in this case.
\end{theorem}
Assumption~\ref{A2.2} in the above theorem can be relaxed provided 
the state space $\cS$ is communicating under every stationary policy
and the cost function $c$ is supported on a finite set. For more details,
see \cite{MR2595908}. Some other works that also study
ERSC problem with bounded cost functions are \cite{MR1629020,MR1795349}.
Since the simultaneous Doeblin condition in Assumption~\ref{A2.2} is
quite restrictive, we are going to impose some structural 
condition on the cost function, known as {\it near-monotonicity}, which
also allows unbounded cost functions.

\begin{definition}\label{D2.2}
We say that the one-step cost function $c$ is near-monotone with
respect to $\rho$ if
$$\liminf_{x\to\infty}\, \min_{u\in\Act(x)} c(x, u)>\rho.$$
\end{definition}
Suppose that for some stationary Markov control $\tilde\zeta$, we have $\sE_x(c, \tilde\zeta)$ independent of $x\in\cS$ and $c$ is near-monotone with respect to $\sE_x(c, \tilde\zeta)$. It is then shown in 
\cite{MR1709324} that, for $\gamma>0$, there exists a positive $\psi:\cS\to (0, \infty]$ satisfying 
\begin{equation}\label{E2.7}
e^{\gamma\lambda^{*,{\rm m}}_{\rm m}}\, \Psi(x)\geq \inf_{u\in\Act(x)} \left[e^{c(x, u)} \sum_{y\in\cS} \Psi(y) P(y|x, u)\right]\quad \text{for all}\; x\in \cS,
\end{equation}
where $\lambda^{*,{\rm m}}_{\rm m}$ is given by 
\begin{equation}\label{E2.8}
\lambda^{*,{\rm m}}_{\rm m}=\inf_{x\in\cS}\, \inf_{\zeta\in\Um}\sE_x(c, \zeta).
\end{equation}
Furthermore, if $\zeta^*$ is a minimizing selector of \eqref{E2.7}, then $\lambda^{*,{\rm m}}_{\rm m}=\sE_x(c, \zeta^*)$ for all $x\in \{\Psi<\infty\}$. The main idea in
\cite{MR1709324} (motivated from \cite{MR1466928}) is to transform the risk-sensitive minimization problem to a risk-neutral game problem using a change of variables (a `logarithmic transformation' that we see later)
and then use the approach of discounted-control problems for the ergodic risk-neutral game to construct a solution for \eqref{E2.7}.

\begin{definition}
We say a function $F:\cS\to \RR$ is {\it norm-like} if for each integer $n$ the set $\{F\leq n\}$ is either empty or finite. 
\end{definition}
Around the same time multiplicative ergodic theorems with norm-like potential functions $F$ are studied in \cite{MR1787128}. The ideas of \cite{MR1787128} are extended to 
study ERSC problems for norm-like cost function $c$ in \cite{MR1886226}. To explain the result of \cite{MR1886226} we introduce some additional notations. Fix a state $z\in\cS$. For
a Markov policy $\zeta\in \Um$, define
\begin{equation}\label{A01}
\Lambda(\zeta)=\inf\left\{\Lambda\, :\, \Exp_z^\zeta\left[e^{\sum_{t=0}^{\uuptau_z-1} \gamma(c(X_t, \zeta_t)-\Lambda)}\leq 1\right]\right\},
\quad \text{and}\quad \Lambda^*= \inf_{\zeta\in\Um}\Lambda(\zeta).
\end{equation}
The first entrance time to the state $z$ is defined as $\sigma_z=\inf\{n\geq 0\, :\, X_n=z\}$. Let us also define, for $x\in\cS$,
\begin{align*}
\Psi_*(x)&\df\inf_{\zeta\in\Um} \Exp_x^\zeta\left[e^{\sum_{t=0}^{\sigma_z}\gamma(c(X_t, \zeta_t)-\Lambda^*)}\right],
\\
w_*(x)&\df\Argmin_{u\in\Act(x)}\, \left(e^{c(x, u)}\sum_{y\in\cS} \Psi_*(y) P(y|x, u)\right).
\end{align*}
The following result is proved in \cite[Theorem~3.6]{MR1886226}
\begin{theorem}\label{T2.5}
Let $\gamma>0$.
Suppose that $\Act(x)$ is finite for all $x$ and $c(\cdot, u)$ is norm-like for all $u\in\Act$. Also assume that the chain $\textbf{X}$ is communicating under every Markov policy and 
aperiodic under any stationary Markov policy. Then, provided $\Lambda^*$ is finite, the following hold.
\begin{itemize}
\item[(i)] $\Lambda^*=\lambda^{*,{\rm m}}_{\rm m}=\sE_x(c, w_*)$ for all $x$ where $\lambda^{*,{\rm m}}_{\rm m}$ is given by \eqref{E2.8}.
\item[(i)] $\Psi_*$ is finite on $\cS$ and 
\begin{equation*}
e^{\gamma\lambda^{*,{\rm m}}_{\rm m}}\, \Psi_*(x)\geq \inf_{u\in\Act(x)} \left[e^{\gamma c(x, u)} \sum_{y\in\cS} \Psi_*(y) P(y|x, u)\right]\quad \text{for all}\; x\in \cS.
\end{equation*}
\end{itemize}
\end{theorem}
The above result requires $\Lambda^*$ to be finite and the chain to be 
aperiodic under each stationary Markov control. Another result of similar
flavor is recently obtained in \cite{MR4429406}, which we state below.
\begin{theorem}\label{T2.6}
In addition to Assumption~\ref{A2.1} let us also assume the following to
hold.
\begin{itemize}
\item[\hypertarget{i}{(i)}] There exists a state $i_0\in\cS$ such that
$$\min_{u\in\Act(i_0)}P(j|i_0, u)>0\quad \text{for all}\; j\neq i_0.$$
\item[(ii)] $\textbf{X}$ is recurrent under each stationary Markov control.
\item[(iii)] $\inf_{\zeta\in\Usm}\sE_x(c, \zeta)<\infty$ for all $x\in\cS$
and $c$ is near-monotone with respect to $\lambda^{*,{\rm m}}_{\rm m}$ in the sense
of Definition~\ref{D2.2}.
\end{itemize}
Then there exists a positive $\Psi$ satisfying
\begin{equation*}
e^{\gamma\lamstrdm}\, \Psi(x)\geq \inf_{u\in\Act(x)} \left[e^{\gamma c(x, u)} \sum_{y\in\cS} \Psi(y) P(y|x, u)\right]\quad \text{for all}\; x\in \cS,
\end{equation*}
and every minimizing selector is an optimal stationary Markov control. Moreover, $\lambda^{*,{\rm m}}_{\rm m}=\lamstrdm$.
\end{theorem}
\cite{MR4429406} also considers the ERSC problem under a blanket stability 
hypothesis but without the near-monotone condition.
\begin{assumption}\label{A2.3}
Let $\textbf{X}$ be irreducible under any stationary Markov control.
In (i) and (ii) below the function $\Lyap$ on $\cS$ takes values in 
$[1, \infty)$ and $\widehat{C}$ is a positive
constant. Assume that one of the following holds.
\begin{itemize}
\item[(i)] For some positive constant $\beta\in (0,1)$ and a finite set 
$C$ it holds that
$$\sup_{u\in\Act(x)}\sum_{y\in\cS} \Lyap(y) P(y|x, u) \leq
(1-\beta)\Lyap(x) + \widehat{C}\Ind_C(x)\quad x\in\cS,$$
and $\gamma\sup_{\sK}c<\uptheta$ where $\uptheta=\log(\frac{1}{1-\beta})$.
\item[(ii)] For a finite set $C$ and a norm-like function $\ell:\cS\to\RR_+$
it holds that 
$$\sup_{u\in\Act(x)}\sum_{y\in\cS} \Lyap(y) P(y|x, u) \leq
(1-\beta(x))\Lyap(x) + \widehat{C}\Ind_C(x)\quad x\in\cS,$$
\end{itemize}
where $1-e^{-\ell(x)}=\beta(x)$. Moreover, the function
$\ell-\gamma\max_{u\in\Act(\cdot)} c(\cdot, u)$ is norm-like.
\end{assumption}
Condition (ii) above is useful for treating ERSC problems with an unbounded 
cost function $c$. The following result is obtained in \cite{MR4429406}.
\begin{theorem}\label{T2.7}
Suppose that Assumption~\ref{A2.1} and ~\ref{A2.3} hold. Also assume the condition $(i)$ of Theorem~\ref{T2.6}. Then we have the following
\begin{itemize}
\item[(i)] There exists a unique, positive $\Psi$, with $\Psi(i_0)=1$,
satisfying
\begin{equation}\label{ET2.7A}
e^{\gamma\lamstrdm}\, \Psi(x)= \inf_{u\in\Act(x)} \left[e^{\gamma c(x, u)} \sum_{y\in\cS} \Psi(y) P(y|x, u)\right]\quad \text{for all}\; x\in \cS.
\end{equation}
\item[(ii)] A stationary Markov control is optimal if and only if it is a
minimizing selector of \eqref{ET2.7A}.
\end{itemize}
\end{theorem}
Theorem~\ref{T2.6} and \ref{T2.7} are proved using a different approach. 
The authors first solve a nonlinear eigenvalue problem on finite sets
containing $i_0$ and then increase the sets to $\cS$. The condition
\hyperlink{i}{(i)} in Theorem~\ref{T2.6} ensures that the limiting eigenfunction $\Psi$ is positive. This condition is recently removed 
in \cite{Wei-Chen} where the authors used the approach of \cite{MR1886226}
(see Theorem~\ref{T2.5})
to define the eigenfunction.

\subsection{General state space}
Next we describe the results known for the general state space. Some of the important works in this direction are
\cite{DiMasi-07,MR1829070,DiMasi-99b,Jaskiewicz-07}. It is natural that one needs to impose additional conditions to ensure
existence of an eigen-pair. We begin by recalling the following result from \cite{MR1829070}.
\begin{theorem}
Let $\cS$ be a complete separable metric space, $\Act(x)=\Act$ for all $x$ and $\gamma>0$. We assume the following to hold.
\begin{itemize}
\item[\hypertarget{A1}{(A1)}] There exists $\delta<1$ such that for all $x, x^\prime\in\cS$, $B\in\sB(\cS)$ and $u, u^\prime\in\Act$ we have
$P(B|x, u)-P(B|x^\prime, u^\prime)\leq \delta$;
\item[\hypertarget{A2}{(A2)}] $\delta e^{\gamma \norm{c}_{\rm sp}}<1$ where $\norm{c}_{\rm sp}$ denotes the span semi-norm of $c$.
\end{itemize}
Then there exists a bounded, positive continuous function $\Psi$ satisfying
\begin{equation}\label{ET2.6A}
e^{\gamma\lamstrdm}\, \Psi(x)=\min_{u\in\Act}\left[e^{\gamma c(x, u)}\int_{\cS} \Psi(y) P(\D{y}|x, u)\right].
\end{equation}
Furthermore, any minimizing selector of the above equation is an optimal stationary Markov control for the ERSC problem, and $\Psi$ is unique, up to a positive
multiplicative constant, in the class $C_b(\cS)$.
\end{theorem}
Note that condition \hyperlink{A2}{(A2)} above requires $\gamma$ to be small. Writing $\varphi=\log\Psi$ we see from above that
$$\gamma\lamstrdm + \varphi(x)= \min_{u\in\Act}\left[\gamma c(x, u)+ \log\int_{\cS} e^{\varphi(y)} P(\D{y}|x, u)\right].$$
Letting 
$$\mathfrak{T}g(x)=\min_{u\in\Act}\left[\gamma c(x, u)+ \log\int_{\cS} e^{g(y)} P(\D{y}|x, u)\right],$$
it is shown in \cite{MR1829070} that   $\{\mathfrak{T}^{n}\textbf{0}, n \geq 0\},$ where $\textbf{0} :=$ the function identically equal to zero,  converges in the space $C_b(\cS)$ with respect to the span semi-norm. The limit of this sequence
gives a fixed point (up to a positive scalar multiplier) which solves \eqref{ET2.6A}. It can be easily checked that by \hyperlink{A1}{(A1)} , $\mathfrak{T}$ is a local contraction \cite[Proposition~2.2]{DiMasi-99b} and therefore,
uniqueness is immediate.
Condition \hyperlink{A2}{(A2)} above was replaced by a more technical condition in \cite{DiMasi-99b} to obtain \eqref{ET2.6A}. Conditions \hyperlink{A1}{(A1)}-\hyperlink{A2}{(A2)} were
replaced by a minorization
condition and certain exponential moment bounds on the hitting time
to a certain compact set in \cite{DiMasi-07} in order to study the optimality
equation \eqref{ET2.6A}.
These results are further extended in \cite{Jaskiewicz-07,MR2585141,Jaskiewicz-07b} to Borel state spaces and for unbounded cost functions. These works study the ERSC control problem using discounted approximation approach which was initiated in \cite{DiMasi-99b}. For $\beta\in(0, 1)$, let $V_\beta$ be a
positive solution to the dynamic programming equation
\begin{equation}\label{E2.10}
V_\beta(x)=\min_{u\in\Act(x)}\left\{e^{\gamma c(x, u)}\int_{\cS} (V_\beta(y))^\beta P(\D{y}|x, u)\right\}\quad x\in\cS.
\end{equation}
$V_\beta$ is basically the discounted value function associated with a
certain dynamic game \cite[Lemma~1]{Jaskiewicz-07}. Under Assumption~\ref{A2.1}, there exists a unique, bounded solution to \eqref{E2.10} whenever $c$ is bounded (cf. \cite[Proposition~4.1]{DiMasi-99b}).
 Under some additional assumptions on 
the transition kernels (cf. \cite[Theorem~4.2]{DiMasi-99b}),
it can be shown that $\frac{V_\beta}{V_\beta(z)}$, $z\in\cS$ is a fixed point, that converges as $\beta\uparrow 1$ to some $\Psi$ satisfying \eqref{ET2.6A} and
$$\gamma^{-1}\lim_{\beta\to 1}(1-\beta)\log V_\beta(z)=\lamstrdm.$$
The above analysis served as the starting point for \cite{Jaskiewicz-07,MR2585141} where the authors allow the cost to be unbounded.
Suppose that
\begin{equation}\label{E2.11}
x\mapsto \Act(x) \; \text{is upper-semicontinuous.}
\end{equation}
Consider $c\geq 0$ and possibly unbounded and $\gamma>0$. 
Then one can solve \eqref{E2.10} for $c_N=\min\{N, c\}$ to obtain
a sequence of $V_{\beta, N}$ for each $\beta\in (0,1)$. Letting $N\to\infty$,
it is then shown that $\lim_{N\to\infty} V_{\beta, N}=V_\beta$, and
(cf. \cite[Lemma~3.1]{MR2585141}, \cite[Lemma~2]{Jaskiewicz-07})
\begin{equation}\label{E2.12}
V_\beta(x)=\min_{u\in\Act(x)}\left\{e^{\gamma c(x, u)}\int_{\cS} (V_\beta(y))^\beta P(\D{y}|x, u)\right\}\quad x\in\cS.
\end{equation}
Define $m_\beta\df \inf_{\cS} V_\beta$. Letting 
$$\tilde{V}_\beta\df \frac{1}{m_\beta} V_\beta,$$
provided $m_\beta>0$,
in \eqref{E2.12} gives
\begin{equation*}
e^{(1-\beta)\log m_\beta}\tilde{V}_\beta(x)=\min_{u\in\Act(x)}\left\{e^{\gamma c(x, u)}\int_{\cS} (\tilde{V}_\beta(y))^\beta P(\D{y}|x, u)\right\}\quad x\in\cS.
\end{equation*}
In \cite{Jaskiewicz-07}, under the assumption \eqref{E2.11} and $\sup_{\beta\in(0,1)} \tilde{V}_\beta<\infty$, it is shown that, for any sequence $\beta_n\to 1$,
$$\lamstrdm=\frac{1}{\gamma}\lim_{n\to \infty} (1-\beta_n)\log m_{\beta_n}, \quad \text{and}\quad \Psi\df\liminf_{n\to\infty} \tilde{V}_\beta, $$
satisfy
$$e^{\gamma\lamstrdm} \Psi(x)\geq \inf_{u\in\Act(x)} \left[e^{c(x, u)} \int_{\cS} \Psi(y) P(\D{y}|x, u)\right]\quad \text{for all}\; x\in \cS,$$
and every minimizing selector is an optimal stationary Markov control. A similar result is also obtained in \cite[Theorem~5.2]{MR2585141} under a milder hypothesis
that requires $m_\beta$ to be finite for all $\beta\in(0, 1)$, but it also assumes (compare it with Definition~\ref{D2.2}) that
$$\{x\in\cS\, :\, \min_{u\in\Act(x)} c(x, u)\leq \lamstrdm + \delta\}$$
to be compact for some $\delta>0$.

\section{Risk-sensitive control of diffusions}
In this section we review some recent progress on the ERSC problem for controlled diffusions. To begin with,  consider the problem for uncontrolled diffusion.
\subsection{Generalized principal eigenvalue}
Let $\textbf{X}=\{X_t\}$ be a diffusion process in $\Rd$ given by
\begin{equation}\label{E3.1}
\D{X}_t= b(X_t) \D{t} + \upsigma(X_t)\D W_t,
\end{equation}
where $b:\Rd\to\Rd$ is the drift vector, $\upsigma:\Rd\to\RR^{d\times d}$ is the diffusion matrix and $W$ is a $d$-dimensional standard Wiener process on a 
complete filtered probability space $(\Omega, \mathfrak{F}, \Prob)$. There exists a unique strong solution of \eqref{E3.1} (see \cite{MR0336813,MR568986,MR1392450}) for every initial data $X_0=x\in\Rd$, $b$ is Borel measurable and $\upsigma$ is locally Lipschitz and locally non-degenerate ,
provided $b,  \upsigma$ have at-most linear growth.
Let $a(x)=\frac{1}{2}\upsigma\upsigma\transp(x)$.
Given a continuous function $c:\Rd\to \RR$ let us define
$$\sE_x(c)=\limsup_{T\to\infty} \, \frac{1}{T}\Exp_x\left[e^{\int_0^T c(X_s)\D{s}}\right].$$
As in Section~\ref{S-DTMC}, the above quantity is related to an eigen-equation which we describe below. We define the extended generator of \eqref{E3.1} as
\begin{equation}\label{E-lin}
\sL f(x) = \trace(a(x) \grad^2 f) + b(x)\cdot\grad f(x).
\end{equation}
\begin{definition}
We say a pair $(\psi, \lambda)\in C^2(\Rd)\times\RR$ is an eigen-pair of $\sL+c$ if $\psi>0$ in $\Rd$ and 
$$\sL\psi(x) + c(x)\psi(x) =\lambda \psi(x) \quad \text{in}\; \Rd.$$
\end{definition}
To understand the relation between $\sE_x$ and an eigen-pair, let us consider the problem in a smooth bounded domain $D$. More precisely, let $\uptau(D)$ be the first exit time from $D$, that is,
$$\uptau(D)=\inf\{t>0\, :\, X_t\notin D\},$$
and we define
$$\sE_x(c, D)\df \limsup_{T\to\infty} \, \frac{1}{T}\Exp_x\left[e^{\int_0^T c(X_s)\D{s}}\Ind_{\{T<\uptau(D)\}}\right] \quad x\in D.$$
It is then well-known that  there exists a $\lambda_D\in\RR$ such that $\sE_x(c, D)=\lambda_D$ for all $x\in D$ \cite{MR425380} and 
for some $\psi_D\in C^2(D)\cap C(\bar{D})$ we have
\begin{equation}\label{E3.2}
\begin{split}
\sL \psi_D + c(x)\psi_D(x) &=\lambda_D \psi_D(x)\quad \text{in}\; D
\\
\psi & >0 \quad \text{in}\; D,
\\
\psi &=0 \quad \text{on}\; \partial D.
\end{split}
\end{equation}
Thus $(\psi_D, \lambda_D)$ forms a Dirichlet eigen-pair for $\sL+c$ in $D$. Furthermore, $\lambda_D$ is the generalized principal
eigenvalue in the sense of \cite{MR1258192,NP92,PW66}, that is, 
\begin{equation}\label{E-prin}
\lambda_D=\inf\{\lambda\, :\, \exists\, \psi\in C^2_+(D)\cap C(\bar{D})\, \text{satisfying} \, \sL\psi + c(x)\psi\leq \lambda\psi\, \text{in}\, D  \},
\end{equation}
where $C^2_+(D)$ denotes the subset of $C^2(D)$ containing functions that are positive in the interior of $D$. (We define $\lambda_D$ likewise for unbounded $D$.) In addition, it can be easily shown that
the principal eigenfunction $\psi_D$ in \eqref{E3.2} is unique up to a multiplicative constant. So we may want to ask whether $\sE_x(c)=\lambda_{\Rd}$. The answer to this question is negative in general. In fact, \cite[Example~3.1]{MR3926044} shows that
$\lambda_{\Rd} < \inf_x \sE_x(c)$. Thus the risk-sensitive problem in the whole space becomes quite delicate. Let us now recall the 
definition of principal eigenvalue in $\Rd$ from \cite{MR3340379}. The generalized principal eigenvalue of $\sL+c$ in $\Rd$ is defined as follows
\begin{equation}\label{A02}
\lambda_{\Rd}=\inf\{\lambda\, :\, \exists\, \text{positive}\, \psi\in C^2(\Rd)\, \text{satisfying} \, \sL\psi + c(x)\psi\leq \lambda\psi\; \text{in}\; \Rd  \}.
\end{equation}
To illustrate explicit dependence on the potential $c$ we would also use the notation $\lambda_{\Rd}(c)$.
Let us also recall the following result from \cite[Theorem~1.4]{MR3340379}.
\begin{theorem}\label{T3.1}
For every $\lambda\in [\lambda_{\Rd}, \infty)$ there exists a positive $\psi\in C^2(\Rd)$ satisfying 
$$\sL\psi(x) + c(x)\psi(x) =\lambda \psi(x) \quad \text{in}\; \Rd.$$
\end{theorem}
In particular, there are infinitely many eigen-pairs for $\sL+c$ in $\Rd$. To equate  $\sE_x(c)$ with the generalized principal eigenvalue $\lambda_{\Rd}$, we must impose
additional conditions on the diffusion coefficients. More discussion in this direction can be found in \cite{MR3926044}. Another important concern is the simplicity of the principal eigenvalue
$\lambda_{\Rd}$. We need following definition for this purpose.
\begin{definition}[Minimal growth at infinity]
An eigen-pair $(\psi, \lambda)$ is said to have a {\it minimal growth at infinity} if for any compact set $K$ and any positive $v\in C^2(K^c)\cap C(\Rd)$ satisfying
$$\sL v + (c-\lambda) v\leq 0 \quad\text{in} \; K^c,$$
we have $v\geq \kappa \psi$ for some $\kappa>0$.
\end{definition}
The above criterion was introduced by Agmon in \cite{Agmon83} and is very useful to establish simplicity of eigenvalues. 
Let $C^+_\circ(\Rd)$ denote the collection of all non-trivial, non-negative, continuous functions which vanish at infinity. The following notions of monotonicity
are introduced in \cite{MR3926044}. 
\begin{definition}
We say the generalized principal eigenvalue $\lambda_\Rd$ is {\it strictly monotone} at $c$ if for some $h\in C^+_\circ(\Rd)$ we have
$\lambda_\Rd(c-h)<\lambda_\Rd(c)$. We say $\lambda_\Rd$ is {\it monotone on the right} at $c$ if for all $h\in C^+_\circ(\Rd)$ we have
$\lambda_\Rd(c)<\lambda_\Rd(c+h)$.
\end{definition}
It is shown in \cite{MR3926044} that strict monotonicity at $c$ implies $\lambda_\Rd(c-h)<\lambda_\Rd(c)$ for all $h\in C^+_\circ(\Rd)$ and therefore, by convexity, it also implies
monotonicity on the right at $c$. The following equivalence criterion is proved in \cite{MR3926044,MR3917226} (see also \cite{ABP2022} for its generalization to weakly-coupled systems).
\begin{theorem}\label{T3.2}
Suppose that $\lambda_\Rd(c)$ is finite. Then the following are equivalent.
\begin{itemize}
\item[(i)] Eigen-pair $(\psi, \lambda_\Rd(c))$ has a minimal growth at infinity.
\item[(ii)] $\lambda_\Rd$ is monotone on the right at $c$.
\item[(iii)] For some compact ball $B$, we have
\begin{equation}\label{ET3.2A}
\psi(x) = \Exp_x\left[e^{\int_0^{\uuptau}(c(X_t)-\lambda_\Rd)\D{t}}\psi(X_{\uuptau})\right] \quad x\in B^c,
\end{equation}
where $\uuptau=\uptau(B^c)$, the first hitting time to $B$.
\end{itemize}
Furthermore, if one of the above holds, then $\lambda_\Rd(c)$ is simple.
\end{theorem}
The analogy between \eqref{E2.5} and \eqref{ET3.2A} should be noted. To characterize the notion of strict monotonicity we need to introduce the {\it twisted diffusion}.
Given an eigen-pair $(\psi, \lambda)$ of $\sL+c$, the twisted diffusion is given by
$$\D{Y}_t = b(Y_t)\D{t} + 2a(Y_t) \grad\log\psi(Y_t) \D{t} + \upsigma(Y_t) \D{W}_t.$$
The twisted process corresponding to a principal eigen-pair is said to be a {\it ground state process} due to its interpretation in physics. The following result can be found in \cite[Theorem~2.1]{MR3926044} 
(see also \cite{MR2206349,MR2837506,MR3095206}).
\begin{theorem}\label{T3.3}
Suppose that $\lambda_\Rd(c)$ is finite. Then
\begin{itemize}
\item[(i)] For every $\lambda>\lambda_\Rd(c)$, the twisted process corresponding to any eigen-pair $(\psi, \lambda)$ is transient.
\item[(ii)] The following are equivalent.
\begin{itemize}
\item[(a)] $\lambda_\Rd$ is strictly monotone at $c$.
\item[(b)] The ground state process is exponentially ergodic.
\end{itemize}
\end{itemize}
\end{theorem}
Let us remark that \cite{MR3926044} requires $c$ to non-negative for Theorem~\ref{T3.2} and \ref{T3.3} to hold, but this restriction on $c$ is removed in \cite{ABP2022}.

\subsection{ERSC for controlled diffusions}
In this section we review ERSC problem for controlled diffusions. We begin with the exponential linear-quadratic model.
\subsubsection{Exponential Linear-quadratic model}\label{S-ELQ}
Exponential linear quadratic model is a risk-sensitive generalization of the classical linear-quadratic model. Such problems are quite central to the optimal investment models, appearing in  mathematical finance
(see \cite{MR1890061,MR2083714,Fleming-00,MR1972534} and references therein).
More precisely, the controlled diffusion is given by (we consider a slightly more general form)
\begin{equation}\label{E3.4}
\D{X}_t = b(X_t)\D{t} + g(X_t, \zeta_t)\D{t} + \upsigma(X_t) \D{W}_t
\end{equation}
where $\zeta_t$ is a progressively measurable process that is {\it non-anticipative} in the sense that for $s<t$,  $W_t-W_s$ is independent of
the completion of the sigma-field generated by $\{X_0, \zeta_r, W_r\, :\, r\leq s\}$ with respect to $(\mathfrak{F}, \Prob)$.
The control process $\zeta$ is generally assumed take values in some Euclidean space $\RR^m$. As before, we denote the set of all
admissible controls by $\Uadm$. Implicitly, we assume that under every admissible control there exists a  unique strong solution to \eqref{E3.4} in the sense that, given $\zeta, W$ as above on a probability space, there exists an a.s.\ unique $X$ satisfying (\ref{E3.4}).
Let $V:\Rd\to [0, \infty)$ and $\phi:\Rd\times\RR^m\to [0, \infty)$ be two given functions. We define
\begin{equation}\label{E3.5}
\lamstrdf=\inf_{x\in\Rd}\limsup_{T\to\infty} \frac{1}{\gamma T} \log J(x, T)\quad \text{where}\quad J(x, T)=\inf_{\zeta\in\Uadm}\Exp_x\left[e^{\gamma\int_0^T (V(X_t)+\phi(X_t, \zeta_t))\D{t}} \right].
\end{equation}
It should be noted that for a given $\gamma>0$, $J(x, T)$ might not be finite for all $T$. This is known as the {\it breakdown} phenomenon. In fact, \cite[Example~1]{MR1784170} shows that
breakdown can actually happen for some large values of $T$. Thus we need to impose conditions on the coefficients to ensure no breakdown 
\cite{MR1485061,MR1372906}. As mentioned before, the above ERSC problem \eqref{E3.5} is related to the nonlinear eigenvalue problem given by
\begin{equation*}
\gamma\lamstrdf \Psi(x) = \trace(a(x)\grad^2\Psi(x)) + b(x)\cdot \grad\Psi(x) + \min_{u\in\RR^m} \{ g(x, u)\cdot \grad \Psi + \gamma\phi(x, z) \Psi(x)\} + \gamma V(x)\Psi(x)
\end{equation*}
for $x\in\Rd$.
Assume $\gamma>0$. Letting $w(x)=\frac{1}{\gamma} \log\Psi(x)$ in the above, we obtain
\begin{equation}\label{E3.6}
\lamstrdf= \trace(a(x)\grad^2 w(x)) + b(x)\cdot \grad\Psi(x) + Q_0(x, \grad w) + V(x),
\end{equation}
where
$$Q_0(x, \xi)= \gamma \xi a(x) \cdot \xi + \min_{u\in\RR^m} \{ g(x, u)\cdot \xi + \phi(x, z) \}.$$
If we choose $a, g, \phi$ in such a way that
\begin{equation}\label{E3.7}
-\upkappa_1\abs{\xi}^2\leq Q_0(x, \xi)\leq -\upkappa_2\abs{\xi}^2\quad x, \xi\in\Rd,
\end{equation}
for some positive constants $\upkappa_1, \upkappa_2$, and 
\begin{equation}\label{E3.8}
\left|\frac{\partial Q_0(x, \xi)}{\partial \xi}\right|\leq \upkappa_3 \abs{\xi}+\upkappa_4, \quad \left|\frac{\partial Q_0(x, \xi)}{\partial x}\right|\leq \upkappa_3 \abs{x}^2+\upkappa_4,
\end{equation}
for some $\upkappa_3, \upkappa_4>0$, we are in the framework of the Hamilton-Jacobi-Isaacs equation of the ergodic type \cite{MR1173120}. More precisely, if $V$ is coercive, the existence and uniqueness of solution to \eqref{E3.6} 
can be obtained from  \cite{MR1173120}. The following result is proved in \cite[Theorem~3.4]{MR1372906}.
\begin{theorem}
We impose the following conditions.
\begin{itemize}
\item[(i)] $\upsigma, b, g, V, \phi$ are smooth and $\upsigma, b$ are Lipschitz. Also, all the derivatives of $\sigma, b, V$ are bounded by $M(1+\abs{x}^k)$ for some $M, k > 0$;
\item[(ii)] $|g(x, z)|\leq \kappa \tilde{g}(z)$ for some locally bounded $\tilde{g}$ and a constant $\kappa$;
\item[(iii)] For some constant $\upkappa_\circ>0$ we have 
$$\xi a(x)\cdot \xi\geq \upkappa_\circ \abs{\xi}^2\quad \xi,x \in\Rd.$$
\item[(iii)] $V$ is coercive and 
$$\lim_{|z|\to\infty}\phi(x, z)=\infty, \quad \lim_{|z|\to\infty}\frac{\abs{g(x, z)}}{\phi(x, z)}=0\quad \text{uniformly in}\, x.$$
\end{itemize}
Also, assume that \eqref{E3.7}-\eqref{E3.8} hold and
$$Q_0(x, \beta \xi)\geq \beta^2 Q_0(x, \xi) - \kappa \beta(1-\beta) \xi a(x)\cdot\xi -\beta(1-\beta) L(x), \quad \beta\in (0, 1),$$
for $\kappa<\upkappa_2$ and some locally bounded function $L$ satisfying
 $$ V(x)-L(x)\to \infty\quad \text{as}\quad |x|\to\infty.$$
Then there exists a unique eigen-pair  $(w, \lambda)\in C^2(\Rd)\times\RR$, $w$ coercive in nature, satisfying \eqref{E3.6}. Furthermore, when
$Q_0(x, \xi)=-\upkappa\, \xi a(x)\xi\transp$ for some $\upkappa>0$,
we have $\lambda=\lamstrdf$, given by \eqref{E3.5}.
\end{theorem}
Similar result can also be found in \cite[Theorem~3.3]{MR1715337} where the authors imposed some structural assumptions on $g$.

\subsubsection{ERSC with a compact action set}
In this section we review the result on ERSC problem when the action set $\Act$ is compact. Let $\textbf{X}=\{X_t\}$ be a controlled diffusion in $\Rd$ 
governed by the It\^{o} equation
\begin{equation}\label{E3.9}
\D{X}_t = b(X_t, \zeta_t)\D{t} + \upsigma(X_t) \D{W}_t
\end{equation}
where $\zeta_t$ is an admissible control in the sense of Section~\ref{S-ELQ}, taking values in a compact metric space $\Act$. We impose the following conditions on the coefficients
to guarantee the existence and uniqueness of solution to \eqref{E3.9}.
\begin{itemize}
\item[\hypertarget{B1}{{(B1)}}] {\it Local Lipschitz continuity}: The functions $b:\Rd\times\Act\to\Rd$ and $\upsigma:\Rd\to\RR^{d\times d}$ are continuous and satisfy
\begin{equation*}
\abs{b(x, u)-b(y, u)} + \norm{\upsigma(x)-\upsigma(y)}\leq C_R |x-y|\quad \forall \; x, y\in B_R,\; \forall\; u\in\Act,
\end{equation*}
for some constant $C_R$, depending on $R>0$, where $B_R$ denotes the ball of radius $R$ centered at $0$.
\item[\hypertarget{B2}{{(B2)}}] {\it Affine growth condition}: There exists a constant $C_0$ such that
$$\max_{u\in\Act}\, [b(x, u)\cdot x]^+ + \norm{\upsigma(x)}^2\leq C_0(1+|x|^2)\quad x\in\Rd.$$

\item[\hypertarget{B3}{{(B3)}}] {\it Local non-degeneracy}: For each $R>0$, there exists $C_R$ satisfying
$$\xi a(x) \cdot \xi\geq C^{-1}_R |\xi|^2\quad \forall\; \xi\in\Rd, \; x\in B_R,$$
where $a(x)=\frac{1}{2}\upsigma\upsigma\transp(x)$.
\end{itemize}
It is well known that under \hyperlink{B1}{(B1)}-\hyperlink{B3}{(B3)}, for any admissible control $\zeta\in\Uadm$ there exists a unique solution of \eqref{E3.9}
\cite[Theorem~2.2.4]{red-book}. As before, a stationary Markov control would correspond to a Borel measurable map from $\Rd$ to $\Act$ and the class of all
stationary Markov controls is denoted by $\Usm$. It is also well known that for every stationary Markov control in $\Usm$ there exists a unique strong solution
to \eqref{E3.9} which is also a strong Markov process \cite{MR0336813,MR568986,MR1392450}. Now consider a continuous function $c:\Rd\times\Act\to [0, \infty)$ which is locally Lipschitz in $x$ uniformly 
with respect to $u\in\Act$. As before, we define the ERSC problem as follows.
\begin{equation}\label{Ergocost}
\lamstrdf= \inf_{x\in\Rd}\,\inf_{\zeta\in\Uadm} \sE_x(c, \zeta) \quad \text{where}\quad \sE_x(c, \zeta) \,\df\, \limsup_{ T\to\infty} \, \frac{1}{\gamma T}\,
\log \Exp_x^{\zeta} \left[e^{\int_{ 0}^{T} \gamma c(X_t, \zeta_t) \D{t}}\right].
\end{equation}
For the remaining part of this section, we discuss the risk-averse problem and therefore, we shall consider $\gamma$ to be positive. As discussed before,
the above ERSC problem corresponds to a nonlinear eigenvalue problem. For this purpose we introduce a family of operators $\sL_u$, parametrized by $u\in\Act$, defined
as follows
$$\sL_u f(x) = \trace(a(x)\grad^2 f(x)) + b(x, u)\cdot\grad f.$$
We shall be interested in an eigenfunction $\Psi\in C^2(\Rd), \Psi>0$, satisfying
\begin{equation}\label{E3.11}
\min_{u\in\Act}\{\sL_u\Psi(x) + \gamma c(x, u)\Psi(x)\}=\gamma\lamstrdf\,\Psi(x)\quad \text{in}\; \Rd.
\end{equation}
The first major contribution for the ERSC of diffusion came from Fleming and McEneaney \cite{Fleming-95a} (see also \cite{Fleming-91}). 
They prove the following in \cite[Theorem~7.2 and ~7.3]{Fleming-95a}.
\begin{theorem}\label{T3.5}
Suppose that $b(\cdot, u), c(\cdot, u)$ are $C^1$ for each $u\in\Act$, $\upsigma$ is constant, and the following hold.
\begin{itemize}
\item[(i)] $c, \nabla_x c$ are bounded. $\gamma>0$.
\item[(ii)] $\grad_x b$ is bounded in $\Rd$.
\item[\hypertarget{iii}{(iii)}] For some $\kappa>0$ we have
$$(x-y)\cdot (b(x, u)-b(y, u))\leq -\kappa \abs{x-y}^2\quad \forall\; x, y\in\Rd, \; u\in\Act.$$
\end{itemize}
There there exists a $\Psi\in C^2(\Rd), \Psi>0,$ satisfying 
$$\min_{u\in\Act}\{\sL_u\Psi(x) + \gamma c(x, u)\Psi(x)\}=\gamma\lamstrdf\,\Psi(x)\quad \text{in}\; \Rd.$$
Furthermore, any measurable selector of the above equation is an optimal stationary Markov control.
\end{theorem}
Apart from the condition \hyperlink{iii}{(iii)} above, the constant diffusion matrix $\upsigma$ also plays a key role in the above result. These two conditions together render 
Lipschitz regularity to $\log\Psi$. More precisely, the authors use a logarithmic transformation to change the risk-sensitive minimization problem to an 
ergodic game problem. Then using the standard method of vanishing discount, they establish the existence of solution to the Hamilton-Jacobi-Isaacs equation for the ergodic game problem.
In order to extend the result to a more general class of $b$ and
$\upsigma$, \cite{MR2174017} considers ERSC problem 
under a periodic setting. This is the content of our next result.
\begin{theorem}
Suppose that $\upsigma=\sqrt{2}I$ and $b, c$ are periodic in 
$x$ variable with period $1$. Also, assume that $b, c$ are Lipschitz in the
$x$ variable. Then there exists a unique, periodic 
$\Psi\in C^2(\Rd), \Psi>0,$ satisfying 
\begin{equation}\label{ET3.6A}
\Delta \Psi(x) + \min_{u\in\Act}\{b(x, u)\cdot\grad\Psi(x) + 
\gamma c(x, u)\Psi(x)\}=\gamma\lamstrdf\,\Psi(x)\quad \text{in}\; \Rd.
\end{equation}
\end{theorem}
The main idea of the proof goes as follows: one starts with an
exponential of the discounted cost defined as (this is actually the continuous version of the approach that appeared in \cite{DiMasi-99b,MR1795349})
$$
w_\alpha(\gamma, x)=\inf_{\zeta\in\Uadm}
\Exp_x\left[\exp\left(\gamma\int_0^\infty e^{-\alpha t}
c(X_t, \zeta_t) \D{t}\right)\right], \quad \alpha\in (0, 1),
$$
and shows that
\begin{equation}\label{E3.13}
-\alpha\gamma\frac{\partial w_\alpha}{\partial \gamma}
+ \Delta w_\alpha + 
\min_{u\in\Act}\{b(x, u)\cdot\grad w_\alpha + \gamma c(x, u) w_\alpha\}=0,
\end{equation}
and $w_\alpha(0, x)=1$. Note that this is a parabolic equation
when we treat $\gamma$ as a variable. Defining 
$u_\alpha=\gamma^{-1}\log w_\alpha$, it is then shown that 
$\alpha u_\alpha, \grad_x u_\alpha$ are globally bounded, 
uniformly in $\alpha$. This helps us to pass the limit in
\eqref{E3.13} to obtain \eqref{ET3.6A}. This idea of \cite{MR2174017}
was then pushed in \cite{MR2727342,MR2679474,MR2679473} to
solve \eqref{ET3.6A} beyond the  periodic setting and under
near-monotone hypothesis.
\begin{definition}
We say that $c$ is {\it near-monotone} with respect to $\rho\in\RR$ 
if it satisfies
$$
\liminf_{|x|\to\infty} \min_{u\in\Act}c(x, u)>\rho.
$$
\end{definition}
In particular, it was proved in \cite{MR2679473} that if $c$ is near-monotone with respect to $\lamstrdf$ and the diffusion \eqref{E3.9} is
recurrent under each stationary Markov control, then there exists a 
positive $\Psi$ satisfying \eqref{ET3.6A} and every measurable selector 
is an optimal stationary Markov control. But the uniqueness 
of $\Psi$ remains an issue.
The approach of \cite{MR2174017,MR2679473} establishes 
$\gamma\lamstrdf$ as an eigenvalue of nonlinear operator in \eqref{ET3.6A}, but in view of Theorem~\ref{T3.1} (for nonlinear operators, see
\cite[Theorem~2.1]{BR2022}), it is hard to
identify $\gamma\lamstrdf$ as the principal eigenvalue. Thus it is important to establish uniqueness of $\Psi$ (up to a positive multiplicative constant)  for the verification 
result of optimal stationary Markov controls. Define the 
nonlinear operator $\cG$ as
\begin{equation}\label{EG}
\cG f (x) = \min_{u\in\Act}(\sL_u f(x) + \gamma c(x, u) f(x)).
\end{equation}
The generalized principal eigenvalue $\lambda_1(\cG)$
 of $\cG$ is defined as before (along the
lines on \cite{MR3340379,BR2022}) :
$$
\lambda_{1}(\cG)=\inf\{\lambda\, :\, \exists\, \text{positive}\, \psi\in C^2(\Rd)\, \text{satisfying} \, \cG\psi\leq \lambda\psi\, \text{in}\, \Rd \}.
$$
A natural question is: under what condition can we show that
$\lambda_1(\cG)=\gamma\lamstrdf$? If we start with the 
Dirichlet generalized eigenvalue problem for $\cG$ on a sequence of
increasing, smooth bounded domains
and let the domains increase to $\Rd$, then applying Harnack's inequality and monotonicity of generalized principal eigenvalues,
it can be shown that the Dirichlet principal eigenvalues converges to 
$\lambda_1(\cG)$. In \cite{MR2808061}, the author applies this idea to
show that $\lambda_1(\cG)\leq \gamma\lamstrdf$, in general, and furthermore,
if $c$ is near-monotone with respect to $\lamstrdf$ and the diffusion \eqref{E3.9} is
recurrent under each stationary Markov control, then there exists a 
positive $\Psi$ satisfying 
$$\cG\Psi=\gamma\lamstrdf\,\Psi\quad \text{in}\; \Rd,$$
 and every minimizing selector 
is an optimal stationary Markov control. Note that the near-monotone criterion penalizes  instability of the process $\textbf{X}$. Thus it is expected that an optimal stationary Markov control would stabilize the process, 
that is, keep it  positive recurrent. Using this fact, the blanket stability hypothesis was removed in \cite{MR3780687}, proving the the following.
\begin{theorem}\label{T3.7}
Assume \hyperlink{B1}{(B1)}-\hyperlink{B3}{(B3)} and also suppose that
$c$ is bounded and is near-monotone with respect to $\lamstrdf$. In addition, suppose that $b, \upsigma$ are bounded, $\upsigma$ is Lipschitz, $a$ is uniformly elliptic and
\begin{equation}\label{ET3.7A}
\max_{u\in\Act}\, \frac{[b(x, u)\cdot x]^+}{|x|}\to 0
\quad \text{as}\; |x|\to\infty.
\end{equation}
Then there exists a
$\Psi\in C^2(\Rd)$ satisfying $\inf_{\Rd}\Psi>0$ and
\begin{equation}\label{ET3.7B}
\min_{u\in\Act}\{\sL_u\Psi(x) + \gamma c(x, u)\Psi(x)\}=\gamma\lamstrdf\,\Psi(x)\quad \text{in}\; \Rd.
\end{equation}
Moreover, the following hold.
\begin{itemize}
\item[(i)] $\lambda_1(\cG)=\lamstrdf=\inf_{\zeta\in\Uadm} \sE_x(c, \zeta)$ for all $x$.
\item[(ii)] If $v\in\Usm$ is a minimizing selector of \eqref{ET3.7B},
then $v$ is stable and is an optimal stationary Markov control.
\item[(iii)] In addition, if we have 
$\lambda^{*,{\rm d}}_{\rm m}(c)<\lambda^{*,{\rm d}}_{\rm m}(c+h)$ for all $h\in C^+_\circ(\Rd)$ where
$$\lambda^{*,{\rm d}}_{\rm m}(c)\df \inf_{x\in\Rd}\,\inf_{\zeta\in\Usm} \sE_x(c, \zeta),$$
then $\Psi$ in \eqref{ET3.7B} is unique up to a positive scalar multiple and
any optimal stationary Markov control is given by a measurable selector
of \eqref{ET3.7B}.
\end{itemize}
\end{theorem}
Note that Theorem~\ref{T3.7} does not impose any stability assumption. Condition \eqref{ET3.7A} is used to show that any minimizing selector is in fact stable. The boundedness assumption on
$b$ and $c$ was relaxed in \cite[Proposition~5.2]{MR4188837} where the authors allowed polynomial growth of $b$ and $c$. The condition of
monotonicity on the right in Theorem~\ref{T3.7}(iii) is not easy to
verify. To tackle this difficulty, an alternative set of conditions has also been used
for the ERSC problems as follows.
\begin{assumption}\label{A3.1}
There exists a positive $\Lyap\in C^2(\Rd)$ with $\inf_{\Rd}\Lyap>0$
such that one of the following holds:
\begin{itemize}
\item[(i)] There exists an inf-compact, positive $\ell\in C(\Rd)$ and
a compact set $\cK$ satisfying
\begin{equation}\label{EA3.1A}
\sup_{u\in\Act}\sL_u \Lyap \leq \bar\kappa \Ind_\cK -\ell \Lyap
\quad \text{in}\; \Rd,
\end{equation}
for some constant $\bar\kappa$, and 
$\ell-\max_{u\in\Act} \gamma c(\cdot, u)$
is inf-compact.
\item[(ii)] For some positive constants $\bar\kappa, \uptheta$
and a compact set $\cK$
we have
\begin{equation}\label{EA3.1B}
\sup_{u\in\Act}\sL_u \Lyap \leq \bar\kappa \Ind_\cK -\uptheta \Lyap
\quad \text{in}\; \Rd,
\end{equation}
and
$$\limsup_{|x|\to\infty}\, \max_{u\in\Act}\gamma c(x, u)<\uptheta.$$
\end{itemize}
\end{assumption}
We remark here that \eqref{EA3.1A} is not possible 
when $a, b$ are bounded \cite[Proposition~2.6]{MR3340379}. This is the reason for introducing \eqref{EA3.1B}. Also, note that Assumption~\ref{A3.1} does not require $c$ to be near-monotone, but imposes
a blanket stability hypothesis on the stationary Markov controls.
Assumption~\ref{A3.1}(ii) was also used in \cite{MR2796095} to
prove the existence of a solution $\Psi$ of \eqref{ET3.7B} and the existence of an optimal stationary Markov control.
Uniqueness and verification results are settled in \cite{MR3926044}
where the authors prove the following.
\begin{theorem}\label{T3.8}
Assume that \hyperlink{B1}{(B1)}-\hyperlink{B3}{(B3)} and 
Assumption~\ref{A3.1} hold. Then there exists a
$\Psi\in C^2(\Rd)$ satisfying $\Psi>0$ and
\begin{equation}\label{ET3.8A}
\min_{u\in\Act}\{\sL_u\Psi(x) + \gamma c(x, u)\Psi(x)\}=\gamma\lamstrdf\,\Psi(x)\quad \text{in}\; \Rd.
\end{equation}
Moreover, the following hold.
\begin{itemize}
\item[(i)] $\lambda_1(\cG)=\lamstrdf=\inf_{\zeta\in\Uadm} \sE_x(c, \zeta)$ for all $x$.
\item[(ii)] If $v\in\Usm$ is a minimizing selector of \eqref{ET3.8A},
then it is an optimal stationary Markov control.
\item[(iii)]$\Psi$ in \eqref{ET3.7B} is unique up to a positive scalar multiple and any optimal stationary Markov control is given by a measurable selector
of \eqref{ET3.8A}.
\end{itemize}
\end{theorem}
Incidentally, the approach of \cite{MR3926044} does not extend to jump diffusions. The eigenvalue approach in \cite{MR3926044,MR3780687} crucially uses
the Harnack inequality to establish the existence of principal eigenfunction of $\cG$ and the Harnack inequality does not hold for the nonlocal equation
with rough kernels (cf. \cite[Example~1.1]{MR3942851}). To tackle this problem, \cite{AB2022} used the Lyapunov function in Assumption~\ref{A3.1} as a barrier
function to bound the Dirichlet principal eigenfunctions. More precisely, under a stability assumption analogous to Assumption~\ref{A3.1}, \cite{AB2022} studies the
ERSC problem for a class for jump diffusions with jump-kernels having finite measure and establishes
a result analogous to Theorem~\ref{T3.8}.

\subsubsection{Connection to $H_\infty$ control}\label{H-infinity}
In this section, we briefly touch upon the connection between $H_\infty$ control and the small noise asymptotics of ERSC problem. Readers are encouraged to
consult the book \cite{MR1353236} to find out more on $H_\infty$ control. Let us start with a (deterministic) nonlinear, controlled dynamical system
$$\D{y}_t= g(y_t, \zeta_t, \xi_t)\D{t}$$
where $\zeta_t$ and $\xi_t$ are two control process, taking values in some subsets $\Act \subset \RR^m$ and $\mathbb{V}\subset\RR^n$, respectively. 
Define 
$$\cA\df L^2_{\rm loc}(\RR_+, \Act), \quad \cB\df L^2_{\rm loc}(\RR_+, \mathbb{V}).$$ 
We choose 
$\zeta$ as a {\it causal} feedback to $\xi$, that is, $\zeta=\alpha(\xi)$ for some $\alpha:\cB\to\cA$ satisfying
$$\text{if for some $t>0$ we have $\xi=\bar\xi$ in $[0, t]$, then $\alpha(\xi)=\alpha(\bar\xi)$ on $[0, t]$.}$$
The class of such causal feedback controls is denoted by $\Uc$.
The $H_\infty$ control problem can be described as follows. Assume that the dynamical system is stable under the control $\zeta\equiv 0$ and
given a {\it response function} $h:\Rd\times\Act\times\mathbb{V}\to\RR_+$, we have a $\upgamma>0$ and a strategy $\alpha\in \Uc$
satisfying, for some starting point $y_0\in\Rd$, 
\begin{equation}\label{E3.19}
\int_0^T h(y_t, \alpha(\xi)(t), \xi_t) \D{t}\leq \upgamma^2 \int_0^T \abs{\xi_t}^2 \D{t}\quad \text{for all}\; T>0, \ \xi\in \cB.
\end{equation}
The least $\upgamma$ satisfying \eqref{E3.19} is called the $H_\infty$ norm. When existence of $\alpha$ is possible, we
say that the $H_\infty$ suboptimal control problem is solvable with disturbance attenuation level $\upgamma$. Note that the above problem  can also be 
studied by considering the value function
$$V_\upgamma(x)=\inf_{\alpha\in\Uc}\, \sup_{\xi\in\cB}\, \sup_{T\geq 0} \int_0^T [h(y_t, \alpha(\xi)(t), \xi_t)-\upgamma^2 \abs{\xi_t}^2 ]\D{t},$$
where $y_0=x$. Note that $V_\upgamma\geq 0$.
The points where $V_\upgamma$ vanishes correspond to the points from where the $H_\infty$ problem is solvable. As shown in \cite{MR1384966}, 
the value function $V_\upgamma$
is a viscosity solution to 
\begin{equation}\label{E3.20}
\sup_{v\in\mathbb{V}}\inf_{u\in\Act}\{g(x, u, v)\cdot \grad V_\upgamma + h(x, u, v)-\upgamma^2 |v|^2 \}=0\quad \text{in}\; \Rd.
\end{equation}
Thus the $H_\infty$ control problem is related to the study of non-negative viscosity solution to \eqref{E3.20}.
In order to understand the connection of \eqref{E3.20} with the ERSC problem, consider the controlled diffusion
\begin{equation}\label{E3.21}
\D{X}_t = b(X_t, \zeta_t)\D{t} + \left(\frac{\varepsilon}{2\upgamma^2}\right)^{\nicefrac{1}{2}}\D{W}_t,
\end{equation}
where $\varepsilon>0$, and $\zeta$ is an admissible control process taking values in $\Act$. Also, letting $\gamma=\varepsilon^{-1}$ in \eqref{Ergocost}, we define
$$\Lambda_\varepsilon= \inf_{x\in\Rd}\,\inf_{\zeta\in\Uadm} \sE_x(c, \zeta),$$
where the controlled diffusion is given by \eqref{E3.21}. As we have seen before, the above ERSC problem corresponds to the 
eigen-equation
$$\varepsilon^{-1} \Lambda_\varepsilon\Psi_\varepsilon(x)= \frac{\varepsilon}{4\upgamma^2}\Delta\Psi_\varepsilon + \min_{u\in\Act}\{b(x, u)\cdot\grad\Psi_\varepsilon + \varepsilon^{-1} c(x, u)\Psi_\varepsilon(x)\}\quad \text{in}\; \Rd.$$
Letting $W_\varepsilon=\varepsilon\log\Psi_\varepsilon$ we obtain
\begin{equation}\label{E3.22}
\Lambda_\varepsilon = \frac{\varepsilon}{2\upgamma^2}\Delta W_\varepsilon + \max_{v\in\Rd}\,\min_{u\in\Act}\left\{(b(x, u)+v)\cdot \grad W_\varepsilon + c(x, u) -\upgamma^2 \abs{v}^2\right\}.
\end{equation}
Thus, if we could show that the family $\{W_\varepsilon\}$ is locally equicontinuous and $\Lambda_\varepsilon\to \Lambda_0$ as $\varepsilon\to 0$ (along some subsequence), then using the stability
of viscosity solutions, it can be shown from \eqref{E3.22} that
\begin{equation}\label{E3.23}
\max_{v\in\Rd}\,\min_{u\in\Act}\left\{(b(x, u)+v)\cdot \grad W_0 + c(x, u) -\upgamma^2 \abs{v}^2\right\}=\Lambda_0,
\end{equation}
where $W_0$ is a limit of $W_\varepsilon$ in the viscosity sense as $\varepsilon\to 0$. If we set $g(x, u, v)=b(x, u) +v$ and $h(x, u, v)=c(x, u)$, then \eqref{E3.23} is same as
\eqref{E3.20} when $\Lambda_0=0$ and $\mathbb{V}=\Rd$. In fact, the following result was proved in \cite[Theorem~2.10]{Fleming-95b} (the control process $\xi$ does not play any role in this result)
\begin{theorem}\label{T3.9}
Suppose that $g(x, u, v)= b(x) + v$ where $b$ satisfies the conditions in Theorem~\ref{T3.5} and $b(0)=0$, $h(x, u, v)=|h_1(x)|^2$ for some $C^1$ function $h_1:\Rd\to\RR^m$ with $h_1(0)=0$ and
$h_1, \partial_{x_i} h_1$ are bounded for all $i=1, 2, \ldots, d$. Then the $H_\infty$ suboptimal control problem is solvable, starting at the point $0$, at the level $\upgamma$, if and only if 
$\lim_{\varepsilon\to 0}\Lambda_\varepsilon=\Lambda_0=0$.
\end{theorem}
The existence of solution to the more general equation \eqref{E3.23} and a  discussion of  $H_\infty$ control can be found in \cite{Fleming-95a,Fleming-01a}, whereas uniqueness
is discussed in \cite{MR1662969}. In the linear-quadratic setting, similar problems are also studied in \cite{MR1453291,MR1612829,MR1715337}. Let us also mention two interesting works
\cite{MR2465707,MR2653896} where \eqref{E3.23} is studied in the framework of max-plus calculus.

\subsection{Generalized Collatz-Wielandt formula}
Consider a non-negative, irreducible matrix $\mathbb{A}\in\RR^{d\times d}$. Then the celebrated Collatz-Wielandt \cite{MR8590,MR35265} formula suggests
\begin{equation}\label{E3.24}
\lambda(\mathbb{A})=\max_{0\leq x=(x_1, \ldots, x_d)}\, \min_{i: x_i>0}\frac{(\mathbb{A}x)_i}{x_i}=
\min_{0\leq x=(x_1, \ldots, x_d)}\, \max_{i: x_i>0}\frac{(\mathbb{A}x)_i}{x_i},
\end{equation}
where $\lambda(\mathbb{A})$ denotes the Perron-Frobenius eigenvalue of $\mathbb{A}$.  An alternate characterization of $\lambda(\mathbb{A})$
can also be given as follows. Write 
$\mathbb{A}=(a_{ij})=DR$ where 
\begin{equation*}
\begin{gathered}
D=\diag[\kappa_1, \ldots, \kappa_d],\quad \kappa_i\df\sum_{j} a_{ij}
\\
R\df(p(j|i)),\quad p(j|i)\df\frac{a_{ij}}{\kappa_i}.
\end{gathered}
\end{equation*}
Let 
$$\sG=\{(\pi, \tilde{P}):\, \text{$\pi$ is the stationary
probability of the stochastic matrix $\tilde{P}=(\tilde{p}(j|i))$}\}.$$
Then the following representation can be found in \cite{MR2571413}
\begin{equation}\label{E-KL}
\log\lambda(\mathbb{A})
=\sup_{(\pi, \tilde{P})\in\sG} \left(\sum_{i} \pi(i)\bigl[\kappa_i-
D_{\rm KL}(\tilde{p}(\cdot|i)||p(\cdot|i))\bigr]\right),
\end{equation}
where $D_{\rm KL}(\cdot||\cdot)$ denotes the Kullback-Leibler divergence defined as
$$
D_{KL}(\tilde{p}(\cdot|i)\|p(\cdot|i,u))=
\begin{cases}
\sum_{j\in\cS}\tilde{p}(j|i)\log\left(\frac{\tilde{p}(j|i)}{p(j|i,u)}\right) & \text{if}
\; \tilde{p}(\cdot|i)\ll p(\cdot|i,u),
\\
\infty & \text{otherwise}.
\end{cases}
$$
 Given the connection between Perron-Frobenius eigenvalue and the 
risk-sensitive limits, it is natural to expect a similar representation for $\lamstrdm$ or $\lamstrdf$. Let us first consider a DTCMC taking values in a
finite set $\cS$ and the action set $\Act$ is also finite. By
$\Pm(\Act)$ we denote the set of all probability vectors on $\Act$.
The cost function $c$ and transition probability can be extended to
$\Pm(\Act)$ in an obvious fashion. In particular, for $v\in \Pm(\Act)$,
we define
$$c(x, v)=\sum_{u\in\Act} c(x, u) v(u),
\quad P(\cdot|i, v)=\sum_{u\in\Act} P(\cdot|i, u)v(u).$$
Also, extend the set $\Usm$ by allowing the controls to take values in
$\Pm(\Act)$. In \cite{9683319}, the authors consider the ERSC problem
(fix $\gamma=1$, for simplicity)
$$\bar\lambda^{*, {\rm m}}=\max_{i\in\cS}\, \inf_{\zeta\in\Usm}\sE_i(c, \zeta).$$
Then a generalization of \eqref{E-KL} is obtained in \cite{9683319} for the controlled problem. In order to state this result, we denote by $\mathcal{Q}$ the 
set of all stochastic matrices $q=(q_{ij})$ satisfying
$$q_{ij}=0\quad \text{if}\;\; \max_{u\in\Act}P(j|i,u)=0.$$
Let $\mathcal{M}_q$ denote the set of all stationary probability
vectors of $q\in\mathcal{Q}$.
\begin{theorem}[\cite{9683319}]\label{T3.10A}
Define
\begin{align}\label{ET3.10AA}
\tilde{c}(i, q, u)&= c(i, u)-D_{\rm KL}(q(\cdot|i)||P(\cdot|i, u))
\quad \text{where}\quad q(j|i)\df q_{ij},\nonumber
\\
\tilde{c}(i, q, v)&=\sum_{u\in\Act} \tilde{c}(i, q, u) v(u),\quad v\in\Pm(\Act),\nonumber
\\
\widehat\Phi(q, v)&=\sup_{\pi\in\mathcal{M}_q}\sum_{i\in\cS}\pi(i)\tilde{c}(i, q, v). 
\end{align}
Then we have
$$\bar\lambda^{*, {\rm m}}=\min_{v\in\Usm}\, \max_{q\in\mathcal{Q}}\widehat\Phi(q, v)=\max_{q\in\mathcal{Q}}\, \min_{v\in\Usm}\widehat\Phi(q,v),$$
and there exists  a saddle point equilibrium point $(q^*, v^*)$ for the above zero-sum game.
\end{theorem}
The readers must have noticed the analogy between \eqref{E-KL} and
\eqref{ET3.10AA}. In fact, \eqref{ET3.10AA} can be seen as an ergodic 
average of $\tilde{c}$ with respect to a 
suitable Markov chain dictated by the
transition probability matrix $q\in\mathcal{Q}$. Thus the control $v$
has no effect on the dynamics but only on the cost function. 
This forms a single controller zero-sum ergodic game \cite{FR}. We revisit this theme later.

For a general state space, Donsker and Varadhan \cite{MR361998} proved the following min-max formula for diffusions.
\begin{theorem}\label{T3.10}
Let $\cX$ be a compact metric space and $\{T_t\}$ be a strongly continuous, positive semigroup on $C(\cX)$ satisfying $T_t 1=1$ for all $t\geq 0$. Let $L$ be the generator
of $T$ and $c$ be any continuous function on $\cX$. Then
\begin{equation}\label{ET3.10A}
\lambda(c)=\inf_{\psi\in \sD^+}\, \sup_{x\in\cX}\frac{L\psi(x) + c(x)\psi(x)}{\psi(x)}=\sup_{\mu\in\Pm(\cX)}\,\inf_{\psi\in\sD^+}\int_{\cX}\frac{L\psi(x) + c(x)\psi(x)}{\psi(x)}\D\mu,
\end{equation}
where $\sD^+$ is the subset of the domain of $L$ containing all positive functions, $\Pm(\cX)$ denotes the set of all Borel probability measures on $\cX$ and
$$\lambda(c)=\lim_{t\to\infty} \frac{1}{T}\log \norm{T^c_t},$$
where $\{T^c_t\}$ denotes the semigroup generated by $L+c$.
\end{theorem}
If we associate the semigroup $\{T_t\}$ with a Markov process $\{X_t\}$ taking values in $\cX$, that is ,
$$T_t f(x) =\Exp_x[f(X_t)],$$
then $\lambda(c)$ is nothing but 
$$\lambda(c)=\lim_{T\to\infty} \frac{1}{T}\log \sup_{x\in\cX}\Exp_x\left[e^{\int_0^T c(X_t)\D{t}}\right].$$
Thus Theorem~\ref{T3.10} gives a Collatz-Wielandt representation to the risk-sensitive value. 

In the context of discrete time Markov chains, the following representation is proved 
in \cite[Theorem~2.2]{MR3629428}. Let\\

\noindent $\cG := \{\eta(dx, du, dy) = \eta_0(dx)\eta_1(du|x)\eta_2(dy|x,u)$ such that $\eta_0(dx)$ is invariant under the transition kernel $\int\eta_2(dy|x,u)\eta_1(du|x)\}$.

\begin{theorem}\label{Ven}
Let $\cS$ be a compact metric space and $\textbf{X}$ be a controlled Markov process with on $\cS$ with a compact metric action space $\Act$ and a continuous transition kernel $(x,u) \in \cS\times\Act \mapsto P(\D{y}|x,u)$. Suppose that Assumption~\ref{A2.1} holds and
the support of $P(\cdot|x)$ is $\cS$ for all $x$. Also, consider a continuous function $c$ on $\cS\times\Act\times\cS$.
 Then there exists a unique $\lambda_1>0$ (the Perron-Frobenius eigenvalue) and a positive $\Psi\in C(\cS)$ such that
$T\Psi=\lambda_1\Psi$ where $T$ is defined as follows.
$$T f(x)=\sup_\varphi\in\Pm(\Act)\int_\cS e^{c(x,u,y)} f(y)\varphi(du) P(\D{y}|x)\quad \text{for}\; f\in C(\cS).$$
Furthermore, the following representations hold for $\lambda_1$
\begin{align*}
\lambda_1 &= \inf_{0<\psi\in C(\cS)}\, \sup_{\mu\in\Pm(\cS)}\, \frac{\int_\cS T\psi(x) \D{\mu}}{\int_{\cS} \psi\D{\mu}}
=\sup_{0<\psi\in C(\cS)}\, \inf_{\mu\in\Pm(\cS)}\, \frac{\int_\cS T\psi(x) \D{\mu}}{\int_{\cS} \psi\D{\mu}} \\
\log \lambda_1&= \sup_{\eta\in\cG}\Bigg(\int\int\int\eta(dx,du,dy)c(x,u,y) - \\
& \ \int\int\eta_0(dx)\eta_1(du|x)D(\eta_2(dy|x,u)\|p(dy|x,u))\Bigg).
\end{align*}
\end{theorem}
The results of \cite{MR3629428} go far beyond the above setting where  the representation is proved for the optimal value corresponding to a risk-reward problem.
A similar representation is also possible for $\lamstrdm$, the optimal value of 
the ERSC problem. In fact, the following min-max formula is established in \cite[Theorem~3.1]{MR3846082} for DTCMC.
\begin{theorem}
Let $\cS$ be a denumerable state space and we consider the setting of Section~\ref{S-DTMC}. Let $c\geq 0$ and $\gamma>0$. Suppose that Assumption~\ref{A2.1} and ~\ref{A2.2} hold, and
every state $x$ is accessible from $z$, under every stationary Markov policy. Then, if $\lamstrdm$ is finite, we have
\begin{equation}\label{ET3.12A}
\lamstrdm=\inf\{\lambda\; :\; \exists \, \text{positive vector}\, \ \psi\, \text{satisfying}\; e^{\gamma\lambda}\psi(x)\geq \min_{u\in\Act(x)} e^{\gamma c(x, u)}\sum_{y\in\cS}\psi(y) P(y|x, u) \}.
\end{equation}
Furthermore, we have
\begin{equation}\label{ET3.12B}
e^{\gamma\lamstrdm}=\inf_{\psi>0}\, \sup_{x\in\cS} \frac{\min_{u\in\Act(x)} e^{\gamma c(x, u)}\sum_{y\in\cS}\psi(y) P(y|x, u)}{\psi(x)}.
\end{equation}
\end{theorem}
To be precise, \cite{MR3846082} established that for some positive vector $\psi_*$ one has
$$e^{\gamma\lamstrdm}=\sup_{x\in\cS} \frac{\min_{u\in\Act(x)} e^{\gamma c(x, u)}\sum_{y\in\cS}\psi_*(y) P(y|x, u)}{\psi_*(x)}
\geq \inf_{\psi>0}\, \sup_{x\in\cS} \frac{\min_{u\in\Act(x)} e^{\gamma c(x, u)}\sum_{y\in\cS}\psi(y) P(y|x, u)}{\psi(x)}.$$
But the above inequality cannot be strict. Otherwise, for some $\varepsilon>0$ and $\psi>0$ we would have
$$e^{\gamma(\lamstrdm-\varepsilon)}>\sup_{x\in\cS} \frac{\min_{u\in\Act(x)} e^{\gamma c(x, u)}\sum_{y\in\cS}\psi(y) P(y|x, u)}{\psi(x)},$$
which will contradict \eqref{ET3.12A}. This gives us \eqref{ET3.12B}.

For controlled diffusions, similar representation was studied in \cite{MR3571250} with the help of the nonlinear Krein-Rutman theorem. To present the result of \cite{MR3571250}, we consider
a bounded domain $D\subset \Rd$ with a $C^3$ boundary. The reflected controlled diffusion on $\bar{D}$ is given by
\begin{equation}\label{E3.28}
\begin{split}
\D{X}_t & = b(X_t, \zeta_t) \D{t} + \upsigma(X_t) \D{W}_t - \delta(X_t)\D\xi_t,
\\
\D\xi_t&= \Ind_{\{X_t\in\partial D\}} \D\xi_t,\quad \xi_0=0,
\end{split}
\end{equation}
where $b , \upsigma$ are as before (see \hyperlink{B1}{(B1)}-\hyperlink{B3}{(B3)}), $\zeta\in\Uadm$, and $\delta:\Rd\to\Rd$ is co-normal, that is,
$\delta(x)=2a(x){\rm n}(x)$ where ${\rm n}(x)$ denote the unit outward normal on $\partial D$. As before, we define the ERSC problem as
\begin{equation}\label{Ergorefl}
\lamstrdf= \inf_{x\in\Rd}\,\inf_{\zeta\in\Uadm} \sE_x(c, \zeta) \quad \text{where}\quad \sE_x(c, \zeta) \,\df\, \limsup_{ T\to\infty} \, \frac{1}{\gamma T}\,
\log \Exp_x^{\zeta} \left[e^{\int_{ 0}^{T} \gamma c(X_t, \zeta_t) \D{t}}\right],
\end{equation}
where $(X_t, \zeta_t)$ satisfies \eqref{E3.28}. Define
$$C^2_{\delta, +}(D)=\{\psi\in C^2(\bar{D})\; :\; \psi\geq 0, \; \grad \psi\cdot\delta =0\; \text{on}\; \partial D\}.$$
Also, recall the operator $\cG$ from \eqref{EG}. The following representation of $\lamstrdf$ can be found in \cite[Theorem~2.1]{MR3571250}.
\begin{theorem}\label{T3.13}
There exists a unique pair $(\rho, \Psi)\in\RR\times C^2_{\delta, +}(D)$ satisfying
$$\cG\Psi=\rho\Psi\quad \text{and}\quad \max_{\bar{Q}}\Psi=1.$$
Moreover, $\rho=\lamstrdf$, given by \eqref{Ergorefl}, and the following hold.
\begin{align*}
\lamstrdf & = \inf_{0<\psi\in C^2_{\delta, +}(D)}\, \sup_{\mu\in\Pm(\bar{D})} \int_{\bar{D}} \frac{\cG \psi}{\psi}\D\mu
\\
&= \sup_{0<\psi\in C^2_{\delta, +}(D)}\, \inf_{\mu\in\Pm(\bar{D})} \int_{\bar{D}} \frac{\cG \psi}{\psi}\D\mu.
\end{align*}
\end{theorem}
Note that the ERSC problem in \eqref{Ergorefl} is related to the Nisio semigroup given by
$$S_t f(x)\df \inf_{\zeta\in\Uadm}\Exp_x\left[e^{\int_0^t c(X_s, \zeta_s)\D{s}} f(X_t)\right] \quad f\in C(\bar{D}).$$
Thus, Theorem~\ref{T3.13} can be seen as a generalization of Theorem~\ref{T3.10} to nonlinear semigroup. 
One can also have a similar Collatz-Wielandt formula for the  generalized Dirichlet principal eigenvalue which is defined
replacing $\sL+c$ by $\cG$ in \eqref{E-prin}.
\begin{theorem}[\cite{MR4048004}]\label{T3.14}
Let $D$ be a bounded, smooth domain and $\lambda_D(\cG)$ denote the generalized Dirichlet principal eigenvalue of $\cG$ in $D$. Then we have
\begin{align*}
\lambda_D(\cG)&= \inf_{\psi\in C^2_+(D)}\, \sup_{\mu\in\Pm(D)} \int_D \frac{\cG\psi}{\psi}\, \D\mu
\\
&= \sup_{\psi\in C^2_+(D)\cap C_0(D)}\, \inf_{\mu\in\Pm(D)} \int_D \frac{\cG\psi}{\psi}\, \D\mu.
\end{align*}
\end{theorem}
As pointed out in \cite[Remark~2.2]{MR4048004}, the set $C^2_+(D)\cap C_0(D)$ in the second equality cannot be extended to $C^2_+(D)$. By Theorem~\ref{T3.7} and ~\ref{T3.8}
we know that $\lamstrdf=\lambda_1(\cG)$, the generalized principal eigenvalue of $\cG$ in $\Rd$. So one might expect an analog of  Theorem~\ref{T3.14} for $\lamstrdf$. It turns out that for a linear operator $\sL$ of the form \eqref{E-lin}, one has (cf.\ \cite{MR4048004})
$$\lambda_{\Rd}(c)=\inf_{0<\psi\in C^2(\Rd)}\,
\sup_{\mu\in\Pm(\Rd)}\int_{\Rd}\left(\frac{\sL \psi}{\psi} +c\right)\D\mu.$$
But by considering $\sL f= f^{\prime\prime}-f^\prime$ in $\RR$ and $c=0$,
it is shown in \cite[Example~2.3]{MR4048004} that 
$$\lambda_{\Rd}(c)< \sup_{0<\psi\in C^2_b(\Rd)}\,
\inf_{\mu\in\Pm(\Rd)}\int_{\Rd}\frac{\sL \psi}{\psi}\D\mu
\leq \sup_{0<\psi\in C^2(\Rd)}\,
\inf_{\mu\in\Pm(\Rd)}\int_{\Rd}\frac{\sL \psi}{\psi}\D\mu.$$
Thus we need an additional condition on the operator in order to obtain a 
full Collatz-Wielandt type formula. In \cite{MR4048004}, the following
condition, which is slightly stronger than
Assumption~\ref{A3.1}, is used in order to obtain a
Collatz-Wielandt type formula.

\begin{assumption}\label{A3.2}
There exists a positive $\Lyap\in C^2(\Rd)$ with $\inf_{\Rd}\Lyap>0$
such that one of the following hold:
\begin{itemize}
\item[(i)] There exists an inf-compact, positive $\ell\in C(\Rd)$ and
a compact set $\cK$ satisfying
\begin{equation*}
\sup_{u\in\Act}\sL_u \Lyap \leq \bar\kappa \Ind_\cK -\ell \Lyap
\quad \text{in}\; \Rd,
\end{equation*}
for some constant $\bar\kappa$, and 
$\upbeta\ell-\max_{u\in\Act} \gamma c(\cdot, u)$
is inf-compact, for some $\upbeta\in (0, 1)$.
\item[(ii)] For some positive constants $\bar\kappa, \uptheta$
and a compact set $\cK$,
we have
\begin{equation*}
\sup_{u\in\Act}\sL_u \Lyap \leq \bar\kappa \Ind_\cK -\uptheta \Lyap
\quad \text{in}\; \Rd,
\end{equation*}
and
$$\limsup_{|x|\to\infty}\, \max_{u\in\Act}\gamma c(x, u)<\uptheta.$$
\end{itemize}
\end{assumption}
By $\sorder(\Lyap)$ we denote the class of functions growing slower than
$\Lyap$, that is, $f\in\sorder(\Lyap)$ if and only if
$$\limsup_{|x|\to\infty} \frac{|f(x)|}{\Lyap(x)}=0.$$

\begin{theorem}[\cite{MR4048004}]
Suppose that \hyperlink{B1}{(B1)}-\hyperlink{B3}{(B3)} and
Assumption~\ref{A3.2} hold. Then we have
\begin{align*}
\lamstrdf=\lambda_1(\cG)&=\inf_{0<\psi\in C^2(\Rd)}\, 
\sup_{\mu\in\Pm(\Rd)}\int_{\Rd}\frac{\cG\psi}{\psi}\, \D\mu
\\
&= \sup_{0<\psi\in C^2(\Rd)\cap\sorder(\Lyap)}\, 
\inf_{\mu\in\Pm(\Rd)}\int_{\Rd}\frac{\cG\psi}{\psi}\, \D\mu.
\end{align*}
\end{theorem}

\section{Risk-sensitive control of  continuous time Markov chains}
In this section we review the recent developments on ERSC for continuous time controlled Markov chains. 
We consider a continuous time controlled Markov chain (CTCMC) $\mathbf{X}=\{X_t$\,,$t\ge 0\}$, on a denumerable state space $\cS$, controlled by the control process $\zeta_t$\,, $t\ge 0$\,, taking values in $\Act$.
As before, $\Act$ is the action space of the controller, which is assumed to be a Borel space with Borel $\sigma$-algebra $\sB(\Act)$. For each $i\in S$, let $\Act(i)$ be the space of all admissible actions of the controller when the system is at state $i$.
 Let $\sK \df \{(i,u) : i\in S, u\in \Act(i)\}$ be the set of all feasible state action pairs.
 As before, we denote by $c:\sK \to \RR_{+}$ the running cost function. The transition rates $q(j|i, u)$, $u\in \Act(i)$\,, $i,j\in S$, satisfy the condition $q(j|i,u)\ge 0$ for all $u\in\Act(i), i, j\in S$ and $j\neq i$. In addition,
we also impose the following:
\begin{assumption}\label{A4.1}
\begin{itemize}
\item[(a)] For each $i\in S$, the admissible action space $\Act(i)$ is a nonempty compact subset of $\Act$\,.
\item[(b)] The model is conservative:
\begin{equation*}
\sum_{j\in S} q(j|i,u) = 0\quad\forall\,\, u\in\Act(i),\,\, i\in S\,.
\end{equation*}
\item[(c)] The model is stable: 
\begin{equation*}
q(i) \,\df\, \sup_{u\in\Act(i)} (- q(i|i,u))\,=\sup_{u\in\Act(i)}\sum_{j\neq i} q(j|i, u)\, <\,  \infty\quad\forall \,\, i\in S\,.
\end{equation*}
\end{itemize}
\end{assumption}
For each $i,j\in S$, $q(j|i,u)$ is a measurable map on $\Act(i)$\,. Let $c : S\times\Act\to\RR_{+}$ be the running cost function. 

Following \cite{K85} (see also \cite{GL19,GZ19,PZ20}) we briefly describe the evolution of the CTCMC. Let $\cS_{\infty}\,\df\, \cS\cup \{i_{\infty}\}$ for an isolated point $i_{\infty}\notin \cS$. Define the canonical sample space $\Omega \,\df\, (\cS\times (0, \infty))^{\infty}\cup \{(i_0, \theta_1, i_1,\dots,\theta_m,i_m,\infty,i_{\infty},\infty,i_{\infty},\dots)\mid \theta_k \neq \infty, i_k \neq i_{\infty}\quad\text{for all}\quad 0\le k \le m,\,\, m\ge 1\}\,,$ with Borel $\sigma$-algebra $\sB(\Omega)$\,. For each sample point $\omega = (i_0, \theta_1, i_1,\ldots,\theta_m,i_m,\ldots)\in\Omega$, we set $T_0(\omega) = 0$, $T_k(\omega) = \theta_1 + \theta_2 +\dots +\theta_k$, and define $T_{\infty}(\omega) = \lim_{k\to\infty}T_k(\omega)$\,. Now we define a controlled process $\{X_t\}_{t\ge 0}$ on $(\Omega, \sB(\Omega))$ by
\begin{equation}\label{ES1}
X_t = \sum_{k\ge 0} \Ind_{\{T_{k} \le t < T_{k+1}\}}i_k + \Ind_{\{t\ge T_{\infty}\}}i_{\infty}\quad \text{for}\,\, t\ge 0\,.
\end{equation} 
From \cref{ES1}, it is clear that for any $m\ge 1$ and $\omega \in\Omega$, $T_m(\omega)$ denotes the $m$-th jump moment of the process $X_t$, $i_m$ is the state of the controlled process on $[T_m, T_{m+1})$ and $\theta_m = T_m - T_{m-1}$ denotes the waiting time between jumps (or, sojourn time) at state $i_{m-1}$\,. Also, we add an isolated point $u_{\infty}\notin \Act$ to $\Act$ and let $\Act_{\infty} = \Act\cup \{u_{\infty}\}$ and $\Act(i_{\infty}) = \{u_{\infty}\}$. We do not want to consider the process beyond the time $T_{\infty}$. Thus we assume that $i_{\infty}$ is an absorbing state, that is, $q(j|i_{\infty}, u_{\infty}) = 0$ for all $j\in \cS$. Also, assume that $c(i_{\infty}, u) = 0$ for all $u\in\Act_{\infty}$\,. Consider a filtration $\{\mathfrak{F}_{t}\}_{t\ge 0}$ where $\mathfrak{F}_{t} \,\df\, \sigma((T_{m} \le s, X_{T_m}\in A) : 0\le s \le t,\,\, m \ge 0, A\subset \cS)$, and let $\Tilde{\mathfrak{F}} \,\df\, \sigma(\mathcal{A}\times \{0\}, \mathcal{B}\times (s,\infty) : \mathcal{A}\in \mathfrak{F}_0, \mathcal{B}\in\mathfrak{F}_{s-})$ be the $\sigma$-algebra of predictable sets in $\Omega\times(0, \infty)$ with respect to $\mathfrak{F}_{t}$, where $\mathfrak{F}_{s-} \,\df\, \vee_{t < s} \mathfrak{F}_{t}$.  Also, define $\tilde{H}_0=\cS$ and $\tilde{H}_m=\cS\times((0, \infty]\times\cS_\infty)^m$ for $m\geq 1$.

An admissible policy $\zeta = \{\zeta_t\}_{t\ge 0}$ is a measurable map from $(\Omega\times (0,\infty), \Tilde{\mathfrak{F}})$ to $(\Act_{\infty}, \sB(\Act_{\infty}))$ satisfying
$\zeta_t(\omega)\in\Act(X_{t-}(\omega))$ for all $\omega\in\Omega$ and $t\ge 0$\,. Let $\Uadm$ be the space of all admissible policies. An admissible policy $\zeta$ is said to be a Markov policy if $\zeta_t(w) = \zeta_{t}(X_{t-}(\omega))$ for all $\omega\in\Omega$ and $t\ge 0$\,. The space of all Markov policies is denoted by $\Um$. If the Markov policy $\zeta$ does not have any explicit time dependence, then it is called a stationary Markov policy and $\Usm$ denotes the space of all stationary Markov policies. For $\zeta\in\Uadm$, define
$$\Lambda^m(\D{y}|i_0, \theta_1, \ldots, i_{m}, t)= \tilde{q}(\D{y}|i_m, \zeta_{t+T_m}(i_0, \theta_1, \ldots, i_{m}))\quad m\geq 1,$$
where $\tilde{q}(\D{y}|i, u)$ denotes the non-negative measure on $\cS\setminus\{i\}$ induced by $q(\cdot|i, u)$.
For each $i\in S$ and $\zeta\in \Uadm$, it is well known that (cf.
\cite{K85,GL19,GZ19,MR3175209}) there exists a unique probability measure $\Prob_{i}^{\zeta}$ on $(\Omega, \sB(\Omega))$ such that $\Prob_{i}^{\zeta}(X_0 = i) = 1$ and
\begin{align*}
\Prob_i^\zeta(A_k\times(\D{t}\times\D{y})) &=\int_{A_k} \Prob_{i}^{\zeta}(\D{h}_k) \Ind_{\{\theta_k<\infty\}} \Lambda^k(\D{y}|h_k, t) e^{-\int_0^t \Lambda^k(\cS|h_k, v)\D{v}} \D{t}, \quad A_k\in \sB(\tilde{H}_k), k\geq 1,
\\
\Prob_i^\zeta(A_k\times(\infty, i_\infty)) &=\int_{A_k} \Prob_{i}^{\zeta}(\D{h}_k)\left\{\Ind_{\{\theta_k=\infty\}}+ \Ind_{\{\theta_k<\infty\}}  e^{-\int_0^\infty \Lambda^k(\cS|h_k, v)\D{v}} \right\}, \quad A_k\in \sB(\tilde{H}_k), k\geq 1.
\end{align*}
 Let $\Exp_i^{\zeta}$ be the corresponding expectation operator. Also, from \cite[pp.13-15]{GHL09}, we know that $\{X_t\}_{t\ge 0}$ is a Markov process under any $\zeta\in\Um$ (in fact, strong Markov). 
Under some policies the process $\{X_t\}_{t\ge 0}$ may be explosive. In order to avoid explosion of the CTCMC, we impose the following condition (see \cite{GL19,GZ19},\cite[Assumption~2.2]{GHL09}).
\begin{assumption}\label{A4.2}
There exist a function $\Tilde{\Lyap}: S \to [1,\infty)$ and constants $C_0\neq 0, C_1>0$ and $C_2\ge 0$ such that
\begin{itemize}
\item[(a)] $\sum_{j\in S} \Tilde{\Lyap}(j)q(j|i,u) \le C_0\Tilde{\Lyap}(i) + C_{2}$ for all $(i, u)\in \sK$\,;
\item[(b)] $q(i) \le C_1 \Tilde{\Lyap}(i)$ for all $i\in S$\,.
\end{itemize} 
\end{assumption}
For the rest of this section, we are going to assume that \cref{A4.2} holds. Note that
\cref{A4.2} holds if $\sup_{i\in S} q(i)<\infty$. In this case we can choose $\tilde{\Lyap}$ to
be a suitable constant.
 From \cite[Theorem~3.1]{GP11} (see also, \cite[Proposition~2.2]{GL19}) it also
follows that, under \cref{A4.2}, $\Prob^\zeta_i (T_\infty=\infty)=1$ for all $i\in S$ and $\zeta\in \Uadm$.
 
We also assume the following for our CTCMC (compare with \cref{A2.1}).
\begin{assumption}\label{A4.3}
\begin{itemize}
\item[(a)] For each $i\in \cS$, the map $u\mapsto c(i,u)$ is continuous on $\Act(i)$.
\item[(b)] For each $i\in \cS$ and bounded measurable function $f:\cS\to\RR$, the map $u\mapsto \Sigma_{j\in \cS}f(j)q(j|i,u)$ is continuous on $\Act(i)$.
\end{itemize}
\end{assumption}

For each  admissible control $\zeta$ the ergodic risk-sensitive cost is given by
\begin{equation}\label{EErgoCcost}
\sE_i(c, \zeta) \,\df\, \limsup_{ T\to\infty} \, \frac{1}{\gamma T}\,
\log \Exp_i^{\zeta} \left[e^{\gamma \int_{0}^{T} c(X_t, \zeta_t)\D t}\right],\quad \gamma>0,
\end{equation} 
where $\textbf{X}$ is the CTCMC corresponding to $\zeta$ with initial state $i$.
As before, our aim is to minimize \cref{EErgoCcost} over all admissible policies in $\Uadm$. A policy $\zeta^{*}\in \Uadm$ is said to be optimal if for all $i\in \cS$ 
$$
\sE_i(c,\zeta^{*}) \, = \, \inf_{i\in \cS}\inf_{\zeta\in \Uadm}\sE_i(c, \zeta)\df\lamstrcm \quad \text{for all}\; i\,.
$$
We also define
\begin{equation}\label{E4.3}
\lambda^{*,{\rm c}}_{\rm m}=\inf_{x\in\cS}\, \inf_{\zeta\in\Um}\sE_x(c, \zeta).
\end{equation}
Recall that a stationary Markov process {\bf X} with rate matrix $Q=[q(j|i)]$ is said to be irreducible if for any $i, j \in \cS, i\neq j,$
there exists distinct $i_1, i_2, \ldots, i_k\in \cS$ satisfying $q(i_1|i)\cdots q(j|i_k)>0$ (cf.
\cite[p.~107]{GHL09}). The following result is proved in \cite[Theorem~3.2]{MR3520987} when $\cS$ is finite.
\begin{theorem}\label{T4.1}
Let $\cS$ be finite and Assumption~\ref{A4.1}, \ref{A4.3} hold. Also, assume that the CTCMC \textbf{X} is irreducible under every stationary Markov control in
$\Usm$. Then there exists a positive vector $\Psi$ satisfying
\begin{equation}\label{ET4.1A}
\gamma\lamstrcm \Psi(i) = \min_{u\in\Act(i)}\left[\sum_{j\in S} \Psi(j)q(j|i,u) + \gamma c(i, u)\Psi(i)\right]\quad\text{for}\,\, i\in \cS\,.
\end{equation}
Moreover, any minimizing selector of \eqref{ET4.1A} is an optimal stationary Markov control.
\end{theorem}
\cite[Theorem~3.2]{MR3520987} actually proves the existence of an eigen-pair $(\Psi, \lambda^{*,{\rm c}}_{\rm m})$ (see \eqref{E4.3}) satisfying \eqref{ET4.1A}. Since 
$\cS$ is finite, applying Dynkin's formula it can be easily seen that $\lambda^{*,{\rm c}}_{\rm m}=\lamstrcm$.

For $\cS$ infinite, the first result concerning the ERSC problem appeared in \cite{MR3230073} where the authors proved the existence of an optimal stationary Markov control.
One should note that the ERSC control problem in \cite{MR3230073} was over the set $\Usm$.
\begin{theorem}[\cite{MR3230073}]\label{T4.2}
Suppose that Assumptions~\ref{A4.1} and ~\ref{A4.3} hold. In addition, we also assume the following.
\begin{itemize}
\item[(i)] $\sup_{i\in\cS} q(i)<\infty$;
\item[(ii)] CTCMC \textbf{X} is irreducible under every stationary Markov control;
\item[(iii)] For some state $z\in\cS$, there exist functions $w_1:\cS\to\RR_+$, $w_2:\cS\to [1, \infty)$, $w_2$ is norm-like, and positive
constants $\uptheta, \kappa$ satisfying
\begin{equation}\label{ET4.2A}
e^{-w_1(i)}\sum_{j\in\cS} e^{w_1(j)} q(j|i, u)\leq -\uptheta w_2(i) + \kappa\Ind_{\{z\}}(i)\quad \text{for all}\; i\in\cS, \,  u\in \Act(i).
\end{equation}
\item[(iv)] $\sup_{i\in\cS, \zeta\in\Usm}\Exp_i^\zeta[\uuptau_z]<\infty$, where $\uuptau_z$ denotes the return time to the state $z$, that is,
$$\uuptau_z=\inf\{t>0\; :\; X_t=z\}.$$
\end{itemize}
Then, if we choose $\gamma$ small enough so that $\gamma \sup_{\sK} c(i, u)<\uptheta$ for $\uptheta$ as in \eqref{ET4.2A}, the map
$$\Usm\ni \zeta\to \sE_i(c, \zeta)\quad (see \eqref{EErgoCcost})$$
has a minimizer in $\Usm$.
\end{theorem}
A variant of Theorem~\ref{T4.2} was obtained by \cite{MR3937056} where the authors replace the stability condition \eqref{ET4.2A}
by a simultaneous Doeblin condition and
also allow the transition rates to be unbounded.
\begin{theorem}[\cite{MR3937056}]\label{T4.3}
Let Assumption~\ref{A4.1}- ~\ref{A4.3} to hold. In addition, we also
assume the following.
\begin{itemize}
\item[(i)] For some state $z\in \cS$ and a control $\tilde\zeta\in\Usm$
we have $\sE_z(c, \zeta)<\infty$;
\item[(ii)] CTCMC \textbf{X} is irreducible under every stationary Markov control;
\item[(iii)] (Simultaneous Doeblin condition) There exists $t_0\in\RR_+$ and $\upalpha\in (0,1)$
so that $\Prob^\zeta_i(\uuptau_z\geq t_0)\leq \upalpha$ for all
$i\in\cS$ and $\zeta\in\Usm$.
\end{itemize}
Then there exists a positive function 
$\Psi:\cS\to(0, \infty)$ satisfying
\begin{equation}\label{ET4.3A}
\gamma\lambda^{*, {\rm c}}_{\rm m}\, \Psi(i) \geq \min_{u\in\Act(i)}\left[\sum_{j\in S} \Psi(j)q(j|i,u) + \gamma c(i, u)\Psi(i)\right]\quad\text{for}\,\, i\in \cS\,,
\end{equation}
and every minimizing selector is an optimal control in $\Usm$.
\end{theorem}
Proof of Theorem~\ref{T4.3} in \cite{MR3937056} is based on an approximation procedure. The ERSC problem is first solved
for a cost function $c$ having compact support and bounded transition
rate functions. In particular, the optimality equation \eqref{ET4.3A}
is obtained for this approximate model. Then, passing to the limit,
one obtains \eqref{ET4.3A} for the original system. In a recent
work \cite{MR4276515}, the authors also study the ERSC problem for 
CTCMC and establish the existence of $\Psi$ solving \eqref{ET4.3A}
with equality. Instead of simultaneous Doeblin condition, 
\cite{MR4276515} imposes a Lyapunov type stability condition of
the form \eqref{ET4.2A}. But it also assumes  $\inf_{\cS}\Psi>0$ (see \cite[Assumption~5.1]{MR4276515}) 
which is a bit restrictive in nature. A thorough study of ERSC problem appears in the recent work \cite{MR4429406}
under the assumption of blanket stability.

\begin{assumption}\label{A4.4}
Assume that the CTCMC \textbf{X} is irreducible under every stationary Markov control in
$\Usm$. In (i) and (ii) below, the function $\Lyap$ on $S$ takes values in $[1, \infty)$ and
$\widehat{C}$ is a positive constant. Assume also that one of the following hold.
\begin{itemize}
\item[(i)] For some positive constant $\uptheta$ and a finite set $C$ it holds that 
\begin{equation*}
\sup_{u\in\Act(i)} \sum_{j\in S} \Lyap(j) q(j|i,u) \le  \widehat{C}\Ind_{C}(i) - \uptheta\Lyap(i) \quad \forall \quad i\in S\,. 
\end{equation*}
Also assume that $\gamma\norm{c}_\infty\df\gamma\sup_{\sK}c(i,u) < \uptheta$.
\item[(ii)] For a finite set $C$ and a norm-like function $\ell:S\to \RR_+$ it
holds that
\begin{equation*}
\sup_{u\in\Act(i)} \sum_{j\in S} \Lyap(j) q(j|i,u) \le \widehat{C} \Ind_{\cK}(i) - \ell(i)\Lyap(i) \quad \forall \quad i\in S\,. 
\end{equation*}
Moreover, the function $\ell(\cdot)-\max_{u\in\Act(\cdot)} c(\cdot, u)$ is norm-like.
\end{itemize}
\end{assumption}
The above conditions should be compared with Assumption~\ref{A3.1}.

\begin{theorem}[\cite{MR4429406}]\label{T4.4}
Let Assumptions~\ref{A4.1}--\ref{A4.4} hold. Also assume that there exists $i_0\in \cS$ such that 
\begin{equation}\label{ET4.4A}
q(j|i_0,u) > 0\quad \text{for all} \; j\neq i_0, \; \text{and}\; u\in\Act(i_0).
\end{equation}
Then the following hold.
\begin{itemize}
\item[(i)] There exists a unique positive function $\Psi$, $\Psi(i_0)=1$, satisfying
\begin{equation*}
\gamma\lamstrcm\,\Psi(i) = \min_{u\in\Act(i)}\left[\sum_{j\in \cS} \Psi(j)q(j|i,u) + \gamma c(i, u)\Psi(i)\right]\quad\text{for}\,\, i\in \cS\,.
\end{equation*}
\item[(ii)] A stationary Markov control $v\in\Usm$ is optimal if and only if it satisfies
\begin{equation*}
\min_{u\in\Act(i)}\left[\sum_{j\in \cS} \Psi(j)q(j|i,u) + \gamma c(i, u)\Psi(i)\right]
= \left[\sum_{j\in \cS} \Psi(j)q(j|i,v(i))+ c(i, v(i)) \Psi(i)\right]
\end{equation*}
for all $i\in \cS$.
\end{itemize}
\end{theorem} 
Conditioned \eqref{ET4.4A} can be relaxed to include different classes of CTCMC (see \cite[Remark~3.2]{MR4429406}).
As discussed before for DTCMP and controlled diffusion processes, the ERSC problem has also been studied for near-monotone cost functions (see Definition~\ref{D2.2}).
A work in this direction appears in \cite{MR3131498} where the authors prove the following

\begin{theorem}\label{T4.5}
Let Assumptions~\ref{A4.1}, ~\ref{A4.3} hold and the transition rates are bounded, that is, 
$$\sup_{\cS} q(i)<\infty.$$
We also assume that following to hold.
\begin{itemize}
\item[(i)] $-q(i|i, u)>0$ for all $(i, u)\in\sK$;
\item[(ii)] One of the following hold.
\begin{itemize}
\item[(a)] $q(j|i, u)>0$ for all $i\neq j$ and $u\in\Act(i)$.
\item[(b)] For each $i\in \cS$, there exists a finite set $C_i$ such that $\min_{u\in\Act(i)} q(j|i, u)>0$ for all $j\in C_i$ and $\sup_{u\in\Act(i)}q(j|i, u)=0$ for all $j\notin C_i$.
\end{itemize}
\item[(iii)] CTCMC \textbf{X} is recurrent under every stationary Markov control.
\item[(iv)] For some $\zeta\in\Usm$, we have $\sE_i(c, \zeta)$ is finite for all $i\in\cS$. In addition, $c$ is near-monotone with respect to $\lambda^{*,{\rm c}}_{\rm m}$.
\end{itemize}
Then there exists a positive function $\Psi$, satisfying
\begin{equation*}
\gamma\lamstrcm\,\Psi(i) \geq \min_{u\in\Act(i)}\left[\sum_{j\in \cS} \Psi(j)q(j|i,u) + \gamma c(i, u)\Psi(i)\right]\quad\text{for}\,\, i\in \cS\,,
\end{equation*}
and every minimizing selector is an optimal stationary Markov control.
\end{theorem}
The conditions of Theorem~\ref{T4.5} has been relaxed substantially in \cite[Theorem~3.2]{MR4429406}. 
\begin{theorem}\label{T4.6}
Let Assumptions~\ref{A4.1}--\ref{A4.3} hold and 
$$\inf_{\zeta\in\Usm}\sE_i(c, \zeta)<\infty\quad \forall \; i\in\cS.$$
We also assume the following
\begin{itemize}
\item[(i)] There exists $i_0\in \cS$ such that 
\begin{equation*}
q(j|i_0,u) > 0\quad \text{for all} \; j\neq i_0, \; \text{and}\; u\in\Act(i_0).
\end{equation*}
\item[(ii)] $c$ is near-monotone with respect to $\lambda^{*,{\rm c}}_{\rm m}$.
\item[\hypertarget{T4.6c}{(iii)}] CTCMC \textbf{X} is recurrent under every stationary Markov control.
\end{itemize}
Then there exists a positive function $\Psi$, satisfying
\begin{equation*}
\gamma\lamstrcm\,\Psi(i) \geq \min_{u\in\Act(i)}\left[\sum_{j\in \cS} \Psi(j)q(j|i,u) + \gamma c(i, u)\Psi(i)\right]\quad\text{for}\,\, i\in \cS\,,
\end{equation*}
and every minimizing selector is an optimal stationary Markov control.
\end{theorem}
Condition \hyperlink{T4.6c}{(iii)} above can be relaxed for a class of CTCMC, allowing the possibility that $\textbf{X}$ could be transient for some $\zeta\in\Usm$. In fact,
the following is proved in \cite[Theorem~3.3]{MR4429406}.

\begin{theorem}\label{T4.7}
Let $\cS=\{1, 2, \ldots\}$ and Assumption~\ref{A4.1}--\ref{A4.3} hold. We also let 
$$\inf_{\zeta\in\Usm}\sE_i(c, \zeta)<\infty\quad \forall \; i\in\cS,$$
and 
\begin{itemize}
\item[(i)] $c$ is near-monotone with respect to $\lamstrcm_{\rm m}$. CTCMC \textbf{X} is irreducible under every stationary Markov control.
\item[(ii)] There exists a function $W:S\to [1, \infty)$ satisfying $W(i)\geq i$ for all large $i$ and 
\begin{equation*}
\sup_{u\in\Act(i)}\sum_{j\in S} W(i) q(j|i, u)\leq g(i)\quad \text{for}\; i\in \cS,
\end{equation*}
for some function $g:\cS\to \RR$ satisfying $\lim_{i\to\infty} g(i)=0$. Furthermore, for some
$\eta>0$ we have\footnote{The following appears incorrectly in \cite{MR4429406}}
\begin{equation*}
\min_{u\in\Act(i)}\frac{q(i-1|i, u)}{-q(i|i, u)}\,\geq\, \eta\quad \text{for all}\; i\in \cS.
\end{equation*}

\item[(iii)] $q(\cdot|1, u)$ supported on a finite set $C$, independent of $u\in\Act(1)$. 
For $\sD_n\df\{1, \ldots, n\}$, $v\in\Usm$ and any $j\in\sD_n\setminus\{1\}$ there exists
distinct $i_1, i_2, \ldots, i_k\in \sD_n$ we have
\begin{equation*}
q(i_1|1, v(1))q(i_2|i_1, v(i_1))\cdots q(j|i_k, v(i_k))>0\,.
\end{equation*}
Then there exists a positive function $\Psi$ satisfying
\begin{equation}\label{ET4.7A}
\gamma\lamstrcm\, \Psi(i) \geq \min_{u\in\Act(i)}\left[\sum_{j\in \cS} \Psi(j)q(j|i,u)
+ \gamma c(i, u)\Psi(i)\right]\quad\text{for}\,\, i\in \cS\,.
\end{equation}
Furthermore,  any measurable selector of \eqref{ET4.7A} is an
optimal stationary Markov control.
\end{itemize}
\end{theorem}
Theorems~\ref{T4.4}, ~\ref{T4.6} and ~\ref{T4.7} are proved using an approach similar to Theorem~\ref{T3.8}. More precisely, the ERSC problems are first studied in
bounded domains and then it is shown that the optimality equation has a limit as we increase the domains to $\cS$.

\section{Risk-sensitive maximization problems and beyond}

In this section, we briefly review some other types of optimization problems involving risk-sensitive cost criterion. We mainly discuss the maximization problem and the game 
associated to the risk-sensitive cost.

\subsection{Risk-reward problems}
The readers must have noticed that the risk-sensitive minimization problems are not equivalent to the maximization problems. So it naturally becomes interesting to study them separately.
Surprisingly, work on risk-sensitive reward/maximization has been relatively uncommon.
For DTCMP with a finite state space, the maximization problems are covered by Theorem~\ref{T2.2}. In fact, this corresponds to the case $\gamma<0$. Another recent work to deal the risk-reward 
problem is \cite{MR3629428} where the authors consider the maximization problem for DTCMP with a compact state space $\cS$. More precisely, given a continuous
{\it one-stage reward} function $r:\cS\times\Act\times\cS\to \RR$, the following maximization problem is considered in \cite{MR3629428}.
\begin{equation}\label{risk-reward}
\upbeta^{*,{\rm m}}\df \sup_{x\in\cS}\, \sup_{\zeta\in\Uadm}\sE_x(r, \zeta)\quad \text{where}\quad 
\sE_x(r, \zeta)=\limsup_{T\to\infty}\frac{1}{T}\log\Exp_x^\zeta\left[e^{\sum_{t=0}^{T-1} r(X_t, \zeta_t, X_{t+1})}\right].
\end{equation}
The result in \cite[Theorem~2.2 and 2.3]{MR3629428} is stated below.
\begin{theorem}
Let $\cS$ be a compact metric space.
Suppose that $P(\cdot|x, u)$ has full support for every $(x, u)\in\cS\times\Act$. Then $e^{\upbeta^{*,{\rm m}}}$ is the Perron-Frobenius eigenvalue and there exists some positive $\Psi\in C(\cS)$ satisfying
$$e^{\upbeta^{*,{\rm m}}} \Psi(x)=\sup_{v\in\Pm(\Act)}\int_{\cS}\int_{\Act} P(\D{y}|x, u) v(\D{u}) \Psi(y) e^{r(x, u, y)}\df \cT\Psi(x)\quad x\in\cS.$$
Moreover, there exists an optimal $\zeta\in\Usm$ for the maximization problem \eqref{risk-reward}, and the following Collatz-Wielandt representation hold.
\begin{align*}
e^{\upbeta^{*,{\rm m}}} & = \inf_{0<\psi\in C(\cS)}\, \sup_{\mu\in\Pm(\cS)} \frac{\int \cT \psi \D\mu}{\int \psi \D\mu}
\\
& = \sup_{0<\psi\in C(\cS)}\, \inf_{\mu\in\Pm(\cS)} \frac{\int \cT \psi \D\mu}{\int \psi \D\mu}.
\end{align*}
\end{theorem}
Recall from Theorem \ref{Ven} ( \cite[Theorem~3.2]{MR3629428}) that a variational representation of $\upbeta^*$ similar to Theorem~\ref{T3.10A} (and \eqref{E-KL}) is obtained as a consequence of the above result. It is then extended to the case where $P( \cdot |x,u)$ need not have full support, by using an approximation argument.  
For controlled reflected diffusions in a bounded $C^3$ domain (see \eqref{E3.28}), risk-reward problems were studied in \cite[Theorem~2.1]{MR4188837}. 
In case of $\Rd$, as it turns out, the analysis of risk-reward problems does not differ much under the blanket stability hypothesis (see Assumption~\ref{A4.4}). Furthermore, \cite{MR4188837}
studied the maximization problem under a near-monotone condition which we describe now. 
Consider controlled diffusion as in \eqref{E3.9} and assume \hyperlink{B1}{(B1)}-\hyperlink{B3}{(B3)} to hold. Given a continuous reward function $c:\Rd\times\Act\to\RR$ (not necessarily non-negative),
Lipschitz in the first argument uniformly with respect to the second,
we define the risk-sensitive maximization problem as follows.
\begin{equation}\label{risk-reward-diff}
\upbeta^{*,{\rm d}}\df \sup_{x\in\cS}\, \sup_{\zeta\in\Uadm}\sE_x(c, \zeta)
\quad \text{where}\quad 
\sE_x(c, \zeta)=\limsup_{T\to\infty}\frac{1}{\gamma T}\log\Exp_x^\zeta\left[e^{\int_0^T \gamma c(X_t, \zeta_t)\D{t}}\right], \; \; \gamma>0.
\end{equation}
Let us define
$$\sH f =\max_{u\in\Act}\{\sL_u f + \gamma c(x, u) f(x)\}.$$
For each $n$, it known from \cite{MR2409410} that there exists a unique pair $(w_n, \varrho_n)\in C(\bar{B}_n)\cap C^2(B_n)\times\RR$, $B_n$ being the ball of radius $n$ centered at the origin, satisfying
\begin{equation}\label{E5.2}
\begin{split}
\sH w_n &=\gamma \varrho_n w_n \quad \text{in}\; B_n,
\\
w_n &>0 \quad \text{in}\; B_n
\\
w_n &=0 \quad \text{on}\; \partial B_n, \quad w_n(0)=1.
\end{split}
\end{equation}
$(w_n, \varrho_n)$ is called the  Dirichlet generalized principal eigen-pair of $\sH$ in $B_n$. Moreover, we have $\varrho_n<\varrho_{n+1}$ for all $n\in\NN$. It turns out that $\lim_{n\to\infty}\varrho_n=\lambda_1(\sH)$
where $\lambda_1(\sH)$ denotes the generalized principal eigenvalue 
of $\sH$ in $\Rd$ defined as follows.
$$
\lambda_{1}(\sH)=\inf\{\lambda\, :\, \exists\, \text{positive}\, \psi\in C^2(\Rd)\, \text{satisfying} \, \sH\psi\leq \lambda\psi\, \text{in}\, \Rd \}.
$$
\begin{definition}\label{D5.1}
A continuous reward function $c:\Rd\times\Act\to\RR$, which is bounded from above, is said to be {\it near-monotone for the maximization problem}
if
$$\upbeta^{*,{\rm d}}> \lim_{r\to\infty} \, \sup_{B^c_r\times\Act} c(x,u).$$
\end{definition}
The above definition should be compared with the near-monotonicity 
condition used for the minimization problem (see Theorem~\ref{T3.7}).
We must point out that near-monotonicity criterion in \cite{MR4188837} is defined using the limits of Neumann eigenvalues whereas \cite{AB2022} uses 
$\lambda_1(\sH)$. Also, note that $\upbeta^{*,{\rm d}}$ is bounded from above by
$\sup_{\Rd\times\Act} c$, and therefore finite.

We have the following result from \cite[Theorem~4.1]{MR4284521}
\begin{theorem}\label{T5.2}
Let \hyperlink{B1}{(B1)}-\hyperlink{B3}{(B3)} hold, $c$ is bounded from above and $\abs{b}$ has at most linear growth. Suppose that $c$ is near-monotone for the maximization problem in the
sense of Definition~\ref{D5.1}. Then the following hold.
\begin{itemize}
\item[(i)] $\lambda_1(\sH)=\gamma \upbeta^{*,{\rm d}}$.
\item[(ii)] There exists a unique, bounded, positive $\Phi\in C^2(\Rd)$ satisfying
$$\sH\Phi(x)=\gamma\upbeta^{*,{\rm d}}\,\Phi(x), \quad \
 x\in\Rd, \quad \text{and}\quad \Phi(0)=1.$$
\item[(iii)] A stationary Markov control $v$ is optimal if and only if
$$\max_{u\in\Act}\{b(x, u)\cdot \grad \Phi(x) + \gamma c(x, u)\Phi(x)\}= b(x, v(x))\cdot \grad \Phi(x) + \gamma c(x, v(x))\Phi(x)\quad \text{almost surely in}\; \Rd.$$
\end{itemize}
\end{theorem}

\subsection{Risk-sensitive games}
Let us also mention a few interesting works treating the game problems with ergodic risk-sensitive criterion. For finite state DTCMC, zero-sum games
are studied in \cite{MR3900793} whereas \cite{MR3131320,GGPP22a} consider zero-sum games with a countable state space. Some other works dealing with
zero-sum games include: \cite{MR4058410} for controlled diffusion, \cite{MR4185068} for controlled diffusion restricted to an orthant, \cite{MR4310348} for 
controlled reflected diffusion in a bounded domain, \cite{MR3583765} for DTCMC with a general state space. Non zero-sum games with risk-sensitive 
criterion are considered in \cite{MR4330327,MR3801104,MR4396400,MR4436615,GGPP22b,GKPP22}.

\section{Algorithms}
In this section we review the results on the policy iteration and value iteration for the ERSC problem. Towards this end, we also discuss a few 
recent results on equivalent linear programs and reinforcement learning.

\subsection{Policy iteration}
 Since we already known that under suitable hypotheses, the ERSC problem can have an optimal stationary Markov policy, it is natural to investigate if it can be determined through policy or value iteration techniques.
 The first policy iteration for the ERSC problem was considered in \cite{Howard-71} where the authors used policy iteration algorithm to establish existence of an eigen-pair satisfying \eqref{ET2.2A}.
  To begin with, suppose that the DTCMC $\textbf{X}$ takes values in a finite state space
 $\cS$ and is irreducible under every stationary Markov control. Then the policy iteration algorithm (PIA) can be described as follows.

\begin{algorithm}\label{Alg-6.1}
Policy iteration.
\begin{itemize}
\item[1.] \textit{Initialization:} Set $k=0$ and choose a $\zeta_0\in\Usm$.
\item[2.] \textit{Value determination:} Let $\psi_k$ be the unique positive eigenfunction  satisfying $\psi_k(i_0)=1$ for some prescribed $i_0\in\cS$, for the eigenvalue problem  
$$e^{\gamma \lambda_k}\psi_k= \sum_u\zeta_k(u)\left[e^{\gamma c(x, u)}\sum_{y\in\cS} \psi_k(y) P(y|x, u)\right],\quad x\in\cS,$$
where 
$$\lambda_k=\sE_x(c, \zeta_k).$$
\item[3.] \textit{Policy improvement:} Choose $\zeta_{k+1}\in\Usm$ satisfying
$$\zeta_{k+1}(x)\in \Argmin_{u\in\Act(x)} \left[e^{\gamma c(x, u)}\sum_{y\in\cS} \psi_k(y) P(y|x, u)\right].$$
\end{itemize}
\end{algorithm}

The following result can be found in \cite[Theorem~4.7]{Fleming-97a}
\begin{theorem}
Suppose that $\cS, \Act$ are finite, $\gamma>0$ and $\textbf{X}$ is irreducible under every stationary Markov control. Then Algorithm~\ref{Alg-6.1} converges in finite number of steps, that is,
there exists $m\in\NN$ such that $\{\lambda_k\,: 0\leq k\leq m\}$ forms a strictly deceasing sequence until it reaches $\lambda^*$.
\end{theorem}

When the state space $\cS$ is countably infinite, PIA algorithms are studied under two frameworks: (a) under an assumption of Lyapunov stability, (b) the running cost is near-monotone.
Under a Lyapunov stability condition, PIA was established in \cite{MR4429406}, which was then further improved in \cite{Wei-Chen}.
\begin{theorem}[\cite{MR4429406}]\label{T6.2}
Suppose that Assumption~\ref{A2.1} and ~\ref{A2.3} hold with a norm-like $\Lyap$ and in case of Assumption~\ref{A2.3}(ii), there exists an $\eta\in (0, 1)$ so that $\gamma\max_{u\in\Act(\cdot)}c(\cdot, u)\leq \eta\ell$ in $\cS$.
Also suppose that there exist states $i_0, z_0\in\cS$ satisfying
\begin{equation}\label{ET6.2A}
\inf_{u\in\Act(i_0)} P(j|i_0, u)>0\quad \text{for all}\; j\neq i_0, \quad \text{and}\quad \inf_{u\in\Act(j)} P(z_0|j, u)>0\quad \text{for all}\; j.
\end{equation}
Then Algorithm~\ref{Alg-6.1} converges in the sense that $\{\lambda_k\}$ forms a decreasing sequence and $\lim_{k\to\infty}\lambda_k=\lamstrdm$. Furthermore, 
$\psi_k$ converges to the unique solution $\Psi$ in \eqref{ET2.7A}.
\end{theorem}
The rough idea of the proof of Theorem~\ref{T6.2} goes as follows. Define $c_k(i)=c(i, \zeta_{k+1}(i))$ and
$$\theta_k(i)=1-\frac{1}{\psi_k(i)}\left[e^{\gamma c_{k+1}(i)-\gamma\lambda_k} \sum_{y\in\cS} \psi_k(y) P(y|i, \zeta_{k+1}(i))\right].$$
From Algorithm~\ref{Alg-6.1} it follows that $0\leq\theta_k\leq 1$. Once we have a point-wise bound for the sequence $\{\psi_k\}$, the proof of
Theorem~\ref{T6.2} would follow if we could establish that $\theta_k\to 0$ pointwise, as $k\to\infty$. To attain this goal, it is shown in 
\cite[Theorem~4.1]{MR4429406} that the Markov process $\textbf{Y}^{(k)}$ associated to the twisted kernel $\tilde{P}^{(k)}$, defined as,
$$\tilde{P}^{(k)}(j|i)=\frac{\psi_k(j) P(j|i, \zeta_k(i))}{\sum_{j\in\cS}\psi_k(j) P(j|i, \zeta_k(i))}$$
has a unique stationary probability measure $\pi_k$. Furthermore, $\{\pi_k\}$ is tight and every sub-sequential limit has full support. It is then
show that 
$$\lim_{k\to\infty} \sum_{j\in D} \theta(j) \pi_k(j)=0\quad \text{for all finite sets}\; D.$$
This in turn proves that $\theta_k(i)\to 0$ for all $i$, as $k\to\infty$.
Note that \eqref{ET6.2A} is a bit restrictive. This was used to find a {\it small set} for the Markov chain $\textbf{Y}^{(k)}$. Later \cite{Wei-Chen} followed the same approach 
but without the condition \eqref{ET6.2A} and established the convergence of Algorithm~\ref{Alg-6.1}.

Next we come to the second setting where we do not impose any blanket stability assumption  like Assumption~\ref{A2.3}, but work with the near-monotone structure imposed on the running cost. This is 
done under the setting of Theorem~\ref{T2.5}. To describe the result we need a few more notations. For $\zeta\in\Um$, recall the 
quantity $\Lambda(\zeta)$ from \eqref{A01}. We say $\zeta\in\Um$ is {\it stabilizing} if $\zeta\in\Usm$ and $\Lambda(\zeta)<\infty$.
Fix $z\in\cS$ and define the first hitting time of $z$ as $\sigma_z=\inf\{n\geq 0\, :\, X_n=z\}$. Given a stabilizing policy $\zeta_n$ we define the
{\it relative value function} $h_n$ as follows.
$$h_n(i)=e^{-\gamma(c(z, \zeta_n(z))-\Lambda(\zeta_n))}\Exp_x^{\zeta_n}\left[e^{\sum_{t=0}^{\sigma_z} \gamma(c(X_t, \zeta_n(X_t))-\Lambda(\zeta_n)}\right]\quad i\in\cS.$$
Note that $h_n(z)=1$. A new policy $\zeta_{n+1}$ is then defined through the minimization
$$\zeta_{n+1}(i)\in \Argmin_{u\in\Act(i)} e^{c(x, u)}\sum_{j\in\cS} h_n(j)P(j|i, u)).$$
This generates a sequence of stabilizing policies and relative value functions. We need two more additional assumptions.
\begin{itemize}
\item[\hypertarget{H1}(H1)] There exists a positive $\Psi_*$ satisfying
\begin{equation}\label{EH1}
e^{\gamma\Lambda^*}\Psi_*(x)= \inf_{u\in\Act(x)} \left[e^{\gamma c(x, u)} \sum_{y\in\cS} \Psi_*(y) P(y|x, u)\right]\quad \text{for all}\; x\in \cS,
\end{equation}
where $\Lambda^*$ is given by \eqref{A01}.
\item[\hypertarget{H2}{(H2)}] There exists a minimizing selector $w_*$ of \eqref{EH1} such that the transformed kernel
\begin{equation}\label{EH2}
\breve{P}_*(j|i)\df e^{\gamma(c(i, w_*(i))-\Lambda^*)} \frac{\Psi_*(y)P(y|i, w_*(i))}{\Psi_*(i)}
\end{equation}
is positive recurrent with a unique invariant probability $\breve\pi_*$.
\end{itemize}
We say the DTCMC is {\it skip free} if for each $x\in\cS$ there exists a finite set $N_x$ so that $P(N_x|x, u)=1$ for all $u\in\Act$. Now we are ready to state the result from \cite[Theorem~5.4]{MR1886226}.
\begin{theorem}
Grant the setting of Theorem~\ref{T2.5}. Also, assume that \hyperlink{H1}{(H1)}-\hyperlink{H2}{(H2)} hold,  DTCMC is skip free and the relative value functions satisfy the multiplicative Poisson
equation
$$e^{\gamma\Lambda(\zeta_n)}h_n(x)= e^{\gamma c(x, \zeta_n(x))} \sum_{y\in\cS} h_n(y) P(y|x, \zeta_n(x))\quad \forall \; x\in \cS.$$
Suppose moreover that
\begin{itemize}
\item[(i)] $\breve\pi_*(\bar{h}/\Psi_*)<\infty$ where $\bar{h}(x)=\limsup_{n\to\infty} h_n(x)$;
\item[(ii)] For any limit point $(h_\infty, w_\infty, c_\infty)$ of the sequence $\{(h_n, w_n, c(\cdot, \zeta_n))\, :\, n\geq 0\}$, the multiplicative Poisson equation has a solution
$\tilde{h}_\infty$ with transition kernel $P(j|i, w_\infty(i))$, and the associated transformed kernel $\breve{P}_{\infty}$ (defined in the same manner as in \eqref{EH2})
is positive recurrent with an invariant probability measure $\breve\pi_\infty$ and $\breve\pi_\infty(h_\infty/\tilde{h}_\infty)<\infty$.
\end{itemize}
Then $h_n\to \Psi_*$ and $\Lambda(\zeta_n)\searrow \Lambda^*$, as $n\to\infty$.
\end{theorem}

Let us now define a PIA for the controlled diffusion. For a stationary Markov control $v\in \Usm$ we define the operator
$$\sL_v f(x) = \trace(a(x)\grad^2 f(x)) + b(x, v(x))\cdot\grad f(x) + \gamma c(x, v(x)) f(x),$$
and by $\lambda_1(\sL_v)$ we denote the principal eigenvalue of $\sL_v$ in $\Rd$ (see \eqref{A02}). 
\begin{algorithm}\label{Alg-6.2}
Policy iteration.
\begin{itemize}
\item[1.] \textit{Initialization:} Set $k=0$ and choose a $v_0\in\Usm$;
\item[2.] \textit{Value determination:} Let $\Psi_k$ be a principal eigenfunction in $\mathscr{W}^{2,p}_{\rm loc}(\Rd), p>d,$ satisfying $\Psi_k(0)=1$, and 
$$\sL_{v_k} \Psi_k=\gamma \lambda_k \Psi_k\quad \text{in}\; \Rd,$$
where $\lambda_k=\gamma^{-1}\lambda_1(\sL_{v_k})$;
\item[3.] \textit{Policy improvement:} Choose $v_{k+1}\in\Usm$ satisfying
$$v_{k+1}(x)\in \Argmin_{u\in\Act(x)} \left[b(x, u)\cdot\grad \Psi_k + \gamma c(x, u)\Psi_k\right].$$
\end{itemize}
\end{algorithm}
The following convergence result can be found in \cite[Theorem~3.2]{MR4284521}.
\begin{theorem}
Assume the setting of Theorem~\ref{T3.8} and let $b$ have at most linear growth. Also, assume in case of Assumption~\ref{A3.1}(ii) that $\sup_{\Rd\times\Act}(\gamma c)<\uptheta$.
Then $\lambda_k=\sE_x(c, v_k)$ for all $x$ and the Algorithm~\ref{Alg-6.2} converges, that is, $\lambda_k\searrow \lamstrdf$ (given by \eqref{Ergocost}) and
$\Psi_k$ converges weakly in $\mathscr{W}^{2,p}_{\rm loc}(\Rd), p>d,$ to the unique solution $\Psi$ in \eqref{ET3.8A}.
\end{theorem}
An analogous algorithm for the maximization problem has also been proved in \cite[Theorem~4.2]{MR4284521}.
\begin{theorem}
Assume the setting of Theorem~\ref{T5.2} and let $v_0\in\Usm$ and $\lambda_0=\gamma^{-1}\lambda_1(\sL_{v_0})$ be such that
$$\lambda_0> \lim_{r\to\infty} \, \sup_{B^c_r\times\Act} c(x,u).$$
Generate a sequence of  $\lambda_k$ and $v_k$ as follows. Let $\widehat{\Psi}_k$ be the unique principal eigenfunction 
satisfying 
$$\sL_{v_k}\widehat\Psi_k = \gamma \lambda_k \widehat\Psi_k\quad \text{in}\; \Rd, \quad \widehat\Psi_k(0)=1,$$
where $\lambda_1(\sL_{v_k})=\gamma\lambda_k$. Define
$$v_{k+1}(x)\in \Argmax_{u\in\Act(x)} \left[b(x, u)\cdot\grad \widehat\Psi_k + \gamma c(x, u)\widehat\Psi_k\right].$$
Then $\lambda_k\nearrow \upbeta^{*,{\rm d}}$, defined by \eqref{risk-reward-diff}, and $\widehat\Psi_k$ converges weakly to $\Phi$ in Theorem~\ref{T5.2}.
\end{theorem}

For CTCMC, PIA is studied in \cite{MR3230073,MR4276515} under the assumption that both $\cS$ and $\Act$ are finite sets. We consider the setting of Theorem~\ref{T4.1}.
Also, assume that $\sup_{i\in\cS, \zeta\in\Usm}\Exp_i^\zeta[\uuptau_z]<\infty$, where $\uuptau_z$ denotes the return time to a prescribed state $z$, that is,
$$\uuptau_z=\inf\{t>0\; :\; X_t=z\}.$$
For every $\zeta_k\in\Usm$, we let $\lambda_k=\sE_x(c, \zeta_k)$ (which would be independent of $x$) and define
$$h_k(x)=\Exp_x^{\zeta_k}\left[e^{\int_0^{\uuptau_z}\gamma(c(X_t, \zeta_k(X_t))-\lambda_k)\D{t}}\right]\quad x\in\cS.$$
From \cite{MR3230073} we know that 
\begin{equation}\label{E6.4}
\sum_{j\in\cS} h_k(j) q(j|i, \zeta_k(i)) + \gamma c(i, \zeta_k(i))h_k(i)= \gamma \lambda_k h_k(i) \quad i\in\cS.
\end{equation}
As before, the improved policy $\zeta_{k+1}$ is defined through  minimization, that is,
\begin{equation}\label{E6.5}
\zeta_{k+1}(i)\in\Argmin_{u\in\Act(i)}\{\sum_{y\in\cS} h_k(j) q(j|i, u) + \gamma c(i, u)h_k(i)\}.
\end{equation}
Assuming $\cS$ and $\Act$ to be finite, it is shown in \cite[Theorem~5.1]{MR3230073} (see also \cite[Lemma~6.1]{MR4276515}) that the above iteration converges in
finite number of steps and $\zeta_k$ converges to an optimal stationary Markov control. The same PIA above can be extended to countably infinite state space 
under the setting of Theorem~\ref{T4.4}. In fact, the following is proved in \cite[Theorem~4.3]{MR4429406}.
\begin{theorem}
Assume the setting of Theorem~\ref{T4.4} with a norm-like Lyapunov
function $\Lyap$. In case of Assumption~\ref{A4.4}(ii),
let there exist an $\eta\in (0, 1)$ so that $\gamma\max_{u\in\Act(\cdot)}c(\cdot, u)\leq \eta\ell$ in $\cS$. 
Then the PIA (generated by \eqref{E6.4}-\eqref{E6.5}) starting from any $\zeta_0\in\Usm$ converges.
\end{theorem}

\subsection{Relative value iteration}
In this section we review some of the important contributions on value iteration for the ERSC problems. Value iteration (VI) or relative value iteration (RVI) basically provide a
recursive method to generate a sequence of value functions that converge to the solution of the optimality equation. As a by product of this method we can generate
nearly optimal controls. One of the early works dealing with VI appeared in \cite{BHP99}. The authors of \cite{BHP99} 
studied RVI for finite state DTCMC. For simplicity of notation, we restrict ourselves to $\gamma=1$ in this section. Fix a positive function $V_0:\cS\to \RR$ and define a sequence
$\{V_n\}$ recursively  as follows
\begin{equation}\label{E6.6}
V_n(x)=\min_{u\in\Act(x)}\left[ e^{c(x, u)}\sum_{y\in\cS} V_{n-1}(y) P(y|x, u)\right]\quad x\in\cS.
\end{equation}
Fix a point $z\in\cS$.
The relative value functions $\tilde{V}_n$ is defined as $\tilde{V}_n(x)= \frac{V_n(x)}{V_n(z)}$. Also, define the
$n$-th differential cost function as $\lambda_n(x)=\log V_{n}(x)-\log V_{n-1}(x)$.
Then \eqref{E6.6} can be written as 
\begin{equation}\label{E6.7}
\tilde{V}_n(x)=\min_{u\in\Act(x)}\left[ e^{c(x, u)-\lambda_n(z)}\sum_{y\in\cS} \tilde{V}_{n-1}(y) P(y|x, u)\right]\quad x\in\cS.
\end{equation}
Then the following result is proved in \cite{BHP99}.
\begin{theorem}\label{T6.7}
Let $\cS$ be finite and the DTCMC is irreducible under every stationary Markov control. In addition to Assumption~\ref{A2.1}, also suppose that
\begin{equation}\label{ET6.7A}
P(x|x, u)>0\quad \text{for all}\; (x, u)\in\sK.
\end{equation}
Then $\tilde{V}_n(x)\to\Psi(x)$ and $\lambda_n(x)\to\lamstrdm$ for all $x\in\cS$ where $(\Psi, \lamstrdm)$ are given by Theorem~\ref{T2.2}.
\end{theorem}
Though \eqref{ET6.7A} is restrictive, it plays a key role in the analysis of \cite{BHP99}. In particular, this condition is used to establish a contraction
phenomenon for a span semi-norm which is then used to obtain the convergence result for the RVI sequence. Later in \cite{MR2015911} the condition \eqref{ET6.7A} is removed.
The RVI method in \cite{MR2015911} works under a very general set-up. As shown in \cite{MR2015911}, one could transform the given DTCMC model suitably so that \eqref{ET6.7A} holds.
The key hypotheses used by \cite{MR2015911} is as follows.
\begin{itemize}
\item[\hypertarget{H3}{(H3)}] $\cS$ is finite and Assumption~\ref{A2.1} holds. There exists an eigen-pair $(\psi, \lambda), \psi>0$, satisfying
\begin{equation}\label{E6.9}
e^{\lambda}\psi(x) = \min_{u\in\Act(x)}\left[e^{ c(x,u)}\sum_{y\in\cS} \psi(x)P(y|x,u)\right]\quad\text{for}\,\, x\in \cS.
\end{equation}
\end{itemize}
It is easy to see that $\lambda$ has to be $\lambda^*$. Let us now introduce the transformed model from \cite{MR2015911}. Fix $\alpha\in (0, 1)$ and define
$\breve{c}:\sK\to\RR$ as follows.
$$\breve{c}(x, u)=\log ((1-\alpha) e^{c(x, u)} + \alpha).$$
A transformed transition kernel $\breve{P}$ is also defined as follows
$$\breve{P}(y|x, u)= \frac{(1-\alpha) e^{c(x, u)} P(y|x, u) + \alpha\delta_{xy}}{(1-\alpha) e^{c(x, u)} + \alpha},$$
where $\delta_{xy}$ denotes the Kronecker symbol on $\cS$, that is , $\delta_{xy}=0$ for $x\neq y$ and $\delta_{xx}=1$. Note that
$$\breve{P}(x|x, u)>0\quad \text{for all}\; (x, u)\in\sK,$$
and hence \eqref{ET6.7A} holds for the transformed system. Interestingly, if we let
$$\breve{\lambda}=\log((1-\alpha) e^\lambda + \alpha),$$
then it can be easily checked from \eqref{E6.9} that
\begin{equation}\label{E6.10}
e^{\breve\lambda}\psi(x) = \min_{u\in\Act(x)}\left[e^{\breve{c}(x,u)}\sum_{y\in\cS} \psi(x)\breve{P}(y|x,u)\right]\quad\text{for}\,\, x\in \cS.
\end{equation}
Conversely, if have an eigen-pair $(\psi, \breve\lambda)$ satisfying \eqref{E6.10}, then setting 
$$\lambda=\log\left(\frac{e^{\breve\lambda}-\alpha}{1-\alpha}\right),$$
one recovers an eigen-pair $(\psi, \lambda)$ satisfying \eqref{E6.9}. Therefore, it is natural to investigate RVI for \eqref{E6.10}.
As before, given positive $\breve{V}_0$, we define the sequence $\{(\breve{V}_n, \breve{\lambda}_n)\}$ as follows. Set $W_0=\breve{V}_0$ and let
$$W_n(x) = \min_{u\in\Act(x)}\left[e^{\breve{c}(x,u)}\sum_{y\in\cS} W_{n-1}(x)\breve{P}(y|x,u)\right].$$
For a fixed point $z\in\cS$, we define
$$\breve{V}(x)=\frac{W_n(x)}{W_n(z)}\quad \breve\lambda_n(x)=W_n(x)-W_{n-1}(x).$$
Then the following is proved in \cite[Theorem~4.1]{MR2015911}
\begin{theorem}
Assume \hyperlink{H3}{(H3)}. Then $\breve\lambda_n(x)\to\breve{\lambda}$ as $n\to\infty$, for all $x\in\cS$ where $\breve\lambda$ is given 
by \eqref{E6.10}.
\end{theorem}
A non-stationary version of the above RVI can be found in \cite{MR2203811}. Very recently, RVI is studied in \cite{MR4276001} for DTCMC with a compact state space.
In particular, \cite{MR4276001} assumes the following
\begin{itemize}
\item[\hypertarget{H4}{(H4)}] $\cS$ is a compact Polish space, and for some reference positive measure $\nu$ on $\cS$, with full support, we have
$$\cS\times\Act\ni(x, u)\mapsto P(\D{y}|x, y)=\phi(y|x, y)\nu(\D{y})\in\Pm(\cS).$$
Moreover, $\phi(\cdot|\cdot, \cdot)$ is continuous.
\end{itemize}
From the argument of \cite[Theorem~2.2]{MR3629428} and  \hyperlink{H4}{(H4)}, we can find a positive $\Psi\in C(\cS)$, unique up to a multiplicative positive constant, satisfying
\begin{equation}\label{E6.11}
e^{\lamstrdm}\Psi(x)=\min_{u\in\Act}\left[e^{c(x, u)} \int_\cS \Psi(y) P(\D{y}|x, u)\right]\quad x\in\cS.
\end{equation}
Following the same philosophy as before, we can define the RVI as follows: Let $V_0\in C(\cS)$ be positive. Fix $z\in\cS$ and define
$$V_{n+1}(x)=\frac{1}{V_n(z)}\min_{u\in\Act} \left[e^{c(x, u)} \int_\cS V_{n}(y) P(\D{y}|x, u)\right].$$
\begin{theorem}[\cite{MR4276001}]
Let Assumption~\ref{A2.1} and \hyperlink{H4}{(H4)} hold. Then $V_{n}(x)\to \bar{V}(x) \in C(\cS)$, $\nu$ almost surly, as $n\to\infty$ and $\bar{V}(z)=e^{\lamstrdm}$, where 
$\bar{V}$ satisfies \eqref{E6.11}.
\end{theorem}
The above result extends to controlled diffusion in $\Rd$. We discuss the RVI under Assumption~\ref{A3.2}. \cite{MR4276001} also consider 
value iteration under near-monotone setting, but it also requires some less verifiable conditions. Interested readers may consult \cite[Theorem~3.2]{MR4276001}.
Recall from Theorem~\ref{T3.8} that under \hyperlink{B1}{(B1)}-\hyperlink{B3}{(B3)} and Assumption~\ref{A3.2}, there exists a positive $\Psi\in C^2(\Rd)$ satisfying 
\begin{equation}\label{E6.12}
 \min_{u\in\Act}\{\sL_u\Psi(x) +  c(x, u)\Psi(x)\}=\lamstrdf\,\Psi(x)\quad \text{in}\; \Rd.
\end{equation}
Furthermore, $\Psi$ is unique up to a positive multiplicative constant and every minimizing selector of \eqref{E6.12} is an optimal stationary Markov control. Now we define
$$C^2_{\Lyap, +}(\Rd)\df \{g\in C^2(\Rd)\; :\; g>0\quad \sup_{\Rd}\frac{g}{\Lyap}<\infty\}.$$
It is also known that $\Psi\in C^2_{\Lyap, +}(\Rd)$. Let us now consider the parabolic equation
\begin{equation}\label{E6.13}
\partial_t \bar\Phi(t, x)= \min_{u\in\Act}\{\sL_u \bar\Phi(t, x) + c(x, u)\bar\Phi(t, x)\}-\lamstrdf\,\bar\Phi(t, x),\quad t>0,
\end{equation}
and $\bar\Phi(0, x)=\Phi_0\in C^2_{\Lyap, +}(\Rd)$. If we choose $\Phi_0$ such that for some $\kappa>0$
$$\kappa^{-1}\Psi\leq\Phi_0\leq \kappa \Psi,$$
then from the proof of \cite[Theorem~3.2]{MR4276001}, it can be shown that $\bar\Phi(t, x)\to \kappa\Psi$, for some $\kappa>0$, as $t\to\infty$. But equation \eqref{E6.13} contains $\lamstrdf$ which 
is unknown. To replace $\lamstrdf$, we consider a modified equation as follows.
\begin{equation}\label{E6.14}
\partial_t \Phi(t, x)= \min_{u\in\Act}\{\sL_u \Phi(t, x) + c(x, u)\Phi(t, x)\}-\Phi(t,0)\Phi(t, x)\quad t>0,
\end{equation}
with $\Phi(0, x)=\Phi_0$. As can be easily checked, \eqref{E6.13} and \eqref{E6.14} are related by following relation
\begin{equation}\label{E6.15}
\bar\Phi(t,x)=\Phi(t,x) e^{\int_0^t(\Phi(s, 0)-\lamstrdf)\D{s} }.
\end{equation}
From \eqref{E6.15} one can use $\Phi$ for RVI. Note that 
$$\frac{\bar\Phi(t, x)}{\Phi(t, x)}=\frac{\bar\Phi(t, 0)}{\Phi(t, 0)}.$$
Using \eqref{E6.15} and the above relation,
we obtain
\begin{align*}
\frac{\D}{\D{t}}\frac{\Phi(t, x)}{\bar\Phi(t, x)}
= \lamstrdf - \Phi(t, 0)= \lamstrdf - \bar\Phi(t, 0)\frac{\Phi(t, x)}{\bar\Phi(t, x)}.
\end{align*}
Hence
$$\frac{\Phi(t, x)}{\bar\Phi(t, x)}= e^{-\int_0^t\bar\Phi(s, 0)\D{s}}
+ \lamstrdf\,\int_0^t e^{-\int_\tau^t\bar\Phi(s, 0)\D{s}}\D\tau.$$
Thus, if $\bar\Phi(t, x)$ converges to a positive function, then
$\frac{\Phi(t, x)}{\bar\Phi(t, x)}$ converges to a positive constant,
which in turn proves the convergence of $\Phi(t, x)$. From \eqref{E6.15},
we also get that $\Phi(t, 0)\to\lamstrdf$ as $t\to\infty$.
Therefore, to establish RVI, it is enough to study convergence of
$\bar\Phi$.
In fact, we have the following result from \cite[Theorem~3.4]{MR4276001}.
\begin{theorem}
Let \hyperlink{B1}{(B1)}-\hyperlink{B3}{(B3)} and Assumption~\ref{A3.2} hold.
Suppose that $0<\Phi_0\in C^2(\Rd)$ and $\sup_{\Rd}\left(\Phi_0/\Lyap\right)<\infty$.
Then there exists a constant $\upkappa_0=\upkappa_0(\Phi_0)$ such that the
value iteration $\bar\Phi$ in \eqref{E6.13} converges to $\upkappa_0\Psi$
as $t\to\infty$, uniformly on compact sets.
\end{theorem}

\subsection{Linear programming}
In this section we revisit Theorem~\ref{T3.10A}
and formulate an equivalent linear program using its equivalence to a single controller game.  This is done in \cite{9683319}, building upon the ideas from \cite{BorkarCDC}. (Here we present suitably corrected statements of the results therein.) Consider the
setting of Theorem~\ref{T3.10A}, that is, both $\cS$ and $\Act$ are finite.
Also, fix $\gamma=1$ for simplicity.
We do not assume irreducibility, whence
$\sE_i(c,\zeta)$ can depend on the initial state $i$. Under a stationary policy $\zeta_n = v(X_n)$ for all $ n \geq 0$ for some $v: \cS \mapsto \Act$, $\sE_i(c,\zeta)$ exists as a well defined limit. 
Let 
$$
\bar\lambda^{*,{\rm m}}\df \max_{i\in\cS}\,\min_{\zeta\in\Usm}\sE_i(c,\zeta).
$$ 
Let 
$$
\mathcal{Q}\df\{\text{the set of stochastic matrices $Q = [q(j|i)]_{i,j\in\cS}$ on $\cS$}\}.
$$
Let $\mathcal{Q}(i) := \{q(\cdot|i)\}$, which is a copy of the simplex of probability vectors indexed by $\cS$, for each $i\in\cS$. Also define
$$\tilde{c}(i,q,u) := c(i,u) - D_{KL}(q(\cdot|i)\|P(\cdot|i,u)).$$
Consider a controlled Markov chain 
$\{\tilde{X}_n\}$ on $\cS$ as follows. Its action space at $i\in\cS$ is 
$\mathcal{Q}(i)\times\Act$. The controlled transition probabilities are
$$
\tilde{P}(j|i,q,u) = q(j|i), \ i,j\in\cS.
$$
The running \textit{payoff} is $\tilde{c}(i,u,q)$ as above. We shall consider only stationary policies $v: \cS \mapsto \Act$. Let $\mathcal{M}_q$ denote the set of stationary distributions for $q\in \mathcal{Q}$.
Then the risk-sensitive control problem above is equivalent to a zero sum stochastic game with payoff
$$
\min_{v\in\Usm}\, \max_{q\in\mathcal{Q}}\widehat\Phi(q, v)\quad
\text{where}\quad
\widehat\Phi(q, v)=\sup_{\pi\in\mathcal{M}_q}\sum_{i\in\cS}\pi(i)\tilde{c}(i, q, v).
$$
This is a \textit{single controller game} \cite{Filar}
 in the sense that the transition probabilities are controlled by only one of the controllers, the other controller controls only the payoff. It can be shown that this game has a value and it equals $\bar\lambda^{*,{\rm m}}$ \cite{9683319}.
Then the linear program associated with the ERSC problem  over all stationary policies can be derived as in \cite{V1981} and is given by:

\bigskip

\noindent (LP-P) Minimize $\sum_{i\in\cS}\beta_i$ subject to:

\begin{align*}
\beta_i &\geq \sum_{j\in\cS}q(j|i)\beta_j, \ (i,q)\in\cS\times\mathcal{Q},
\\
V_i &\geq \sum_{u\in\Act}\tilde{c}(i,q,u)y_i(u) - \beta_i + \sum_{j\in\cS}q(j|i)V_j, \ (i,u) \in \cS\times\mathcal{Q},
\\
y_i(u)&\geq 0, \ \sum_{j\in \cS} y_j(u) = 1, \ i \in \cS.
\end{align*}

This is the `primal' linear program. In what follows, we denote $\mu\in \PW :=$ the space of finite non-negative measures on
$$\widetilde{W} := \bigcup_{i\in\cS}(\{i\}\times\mathcal{Q}(i))$$ 
as $\mu(i,dq)$ instead of $\mu(\{i\},dq)$ for notational simplicity. The notation `$\int\cdots d\mu(i,q)$' will indicate integration w.r.t.\ the full measure, whereas `$\int\cdots \mu(i,dq)$' will indicate integration over the second variable with the first variable fixed at $i$. The dual linear program then is:

\bigskip

\noindent (LP-D) Maximize $\sum_{i\in\cS} w_i$ subject to:

\begin{align}
\int_{\widetilde{W}}
(\delta_{ij} - q(j|i))d\mu(i,q) &= 0, \ j\in\cS,  \nonumber 
\\
\int_{\widetilde{W}}(\delta_{ij} - q(j|i))d\nu(i,q) + \mu(j,\mathcal{Q}(j)) &= 1, \ j\in\cS, \nonumber \\
\int_{\mathcal{Q}(j)}\tilde{c}(j,u,q)\mu(j,dq) &\geq w_j, \ j\in\cS, u\in\Act, \label{LPline}\\
j\in\cS, \mu,\nu\in\PW, \nonumber
\end{align}
where $\delta_{ij}$ is the Kronecker delta.

One caveat is that unlike \cite{V1981}, the action spaces here are not both finite - one of them being a probability simplex, is not. Thus one has to go via finitary approximations, using the fact that non-negative measures supported on a dense subset of a Polish space are dense in the space of non-negative measures on that space, with respect to the weak$^*$ topology. See \cite{9683319} for details. 

These linear programs are precisely the counterparts of the primal and dual linear programs for multi-chain average cost Markov decision problems \cite{DF}, \cite{Kal} for this specific single controller game. One can show that these linear programs are feasible and have bounded solutions, and the optimal solution is precisely the value of the two person zero sum game. Furthermore, the optimal stationary policy is optimal among all admissible policies and can be recovered from the dual program. (See \cite{9683319} for details.)

One can reverse engineer the dynamic programming equations from this. These are as follows.

\begin{theorem}
The dynamic programming equation
\begin{align*}
\Psi_i  &= \max_{q\in\mathcal{Q}(i)}\sum_jq(j|i)\Psi_j, 
\\
\Psi_i + \mathcal{V}_i &= \min_{u\in\Act}\max_{q\in B_i}\left[\tilde{c}(i,q,u) + 
\sum_{j\in\cS} q(j|i) \mathcal{V}_j\right], 
\\
\mbox{where} \ & \
  \\
B_i &\df \left\{q \in \mathcal{Q}(i) : \sum_j q(j|i)\Psi_j = \Psi_i\right\},
\end{align*}
has a solution $\{(\Psi_i, \mathcal{V}_i)\}$ with  $\Psi_i = \beta_i \ \forall\  i\in\cS$ where $\{\beta_i\}$ is the solution to (LP-P) and $\beta_i =$ the value of the above zero sum game with initial condition $i$ for all $i$. Furthermore, 
$\bar\lambda^{*,{\rm m}} = \max_i\Psi_i$. 
\end{theorem}

The state space $\cS$ can be partitioned into disjoint subsets as $\cS = \cup_{\ell = 1}^k\cI_\ell$ where $\cI_j := \{j : \Psi_i = \beta_j, \  \forall \ i\in \cI_j\}$, $1 \leq j \leq k$.
Performing the maximization with respect to $q\in\mathcal{Q}(i)$  exactly, one can rewrite the dynamic programming equation as
\begin{align*}
\bar\lambda^{*,{\rm m}} &= \max_i\lambda_i\df \max_i \min_{\zeta\in\Usm}\sE_i(c,\zeta), 
\\
\lambda_i\psi_i &= \min_u\left(\sum_{j\in\cI_i}P(j|i,u)e^{c(i,u)}\psi_j\right), \ i \in \cI_\ell, 1 \leq \ell \leq k, 
\\
\lambda_i &= \min_{B_i^*}\sum_{j\in\cI_i}\left(\frac{P(j|i,u)e^{c(i,u)}\psi_j}{\sum_{j'}P(j'|i,u)e^{c(i,u)}\psi_{j'}}\right)\lambda_{j'}, \ (i,j) \in \cI_\ell, 1 \leq \ell \leq k,
\end{align*}
where ${B}^*_i$ is the set of minimizers in the second equation above. This is the counterpart for risk-sensitive control of the classical result of Howard \cite{Howard} for average cost dynamic programming equation for multi-chain problems. Note the appearance of the so called `twisted kernel' in the last equation. Analogous results are possible  for risk-sensitive reward problems \cite{BorkarCDC} (\cite{9683319} restates them correcting the order of maximization in the dynamic programming equation to as it should be.).

As a spin-off, this allows us to handle risk-sensitive control with risk-sensitive constraints\footnote{This material is new, the details will appear elsewhere.}. Consider minimization of $\sE_i(c, \zeta)$ subject to an additional constraint
\begin{equation}\label{constraint}
\limsup_{T\to\infty}\frac{1}{T}\log\Exp\left[e^{\sum_{t=0}^{T-1}k(X_t,\zeta_t)}\right] \leq C, 
\end{equation}
for a prescribed $k: \cS\times\Act \mapsto \RR$ and a constant $C > 0$. Define
$$\tilde{k}(i,q,u) := k(i,u) - D_{KL}(q(\cdot|i)\|p(\cdot|i,u)).$$  Denote the $\widehat\Phi(q,v)$ above as $\widehat\Phi(q,v,c)$ in order to render explicit its dependence on the per stage cost function $c$. Thus we can also define $\widehat\Phi(q,v,k)$ analogously. This leads to the convex program:\\

\noindent Minimize $\max_{q\in\mathcal{Q}}\widehat\Phi(q,v,c)$ subject to $\max_{q'\in\mathcal{Q}}\widehat\Phi(q',v,k) \leq C.$\\

Then by standard Lagrange multiplier theory \cite{MR0238472}, one can consider the unconstrained minimization of
$$\max_{q\in\mathcal{Q}}\widehat\Phi(q,v,c) + \Gamma(\max_{q'\in\mathcal{Q}}\widehat\Phi(q',v,k) - C),$$
where $\Gamma \geq 0$ is the Lagrange multiplier. The primal program is:\\

\noindent  Minimize $\sum_{i\in\cS}(\beta_i + \Gamma(\beta'_i - C))$ subject to:

\begin{align*}
\beta_i &\geq \sum_{j\in\cS}q(j|i)\beta_j, \ (i,q)\in\cS\times\mathcal{Q},
\\
V_i &\geq \sum_{u\in\Act}\tilde{c}(i,q,u)y_i(u) - \beta_i + \sum_{j\in\cS}q(j|i)V_j, \ (i,q) \in \cS\times\mathcal{Q},
\\
\beta_i' &\geq \sum_{j\in\cS}q'(j|i)\beta_j', \ (i,q')\in\cS\times\mathcal{Q},
\\
V_i' &\geq \sum_{u\in\Act}\tilde{k}(i,q',u)y_i(u) - \beta_i' + \sum_{j\in\cS}q'(j|i)V_j', \ (i,q') \in \cS\times\mathcal{Q},
\\
y_i(u)&\geq 0, \ \sum_j y_j(u) = 1.
\end{align*}

The dual linear program is: \\

\noindent Maximize $\sum_{i,j\in S}w_{i}$ subject to:

\begin{align}
\int_{\widetilde{W}}
(\delta_{ij} - q(j|i))d\mu(i,q) &= 0, \ j\in\cS,\nonumber 
\\
\int_{\widetilde{W}}(\delta_{ij} - q(j|i))d\nu(i,q) + \mu(j,\mathcal{Q}(j)) &= 1, \ j\in\cS, \nonumber \\
\int_{\widetilde{W}}
(\delta_{ij} - q'(j|i))d\mu'(i,q') &= 0, \ j\in\cS,\nonumber 
\\
\int_{\widetilde{W}}(\delta_{ij} - q'(j|i))d\nu'(i,q') + \mu'(j,\mathcal{Q}(j)) &= 1, \ j\in\cS, \nonumber \\
\int_{\mathcal{Q}(i)}\tilde{c}(i,q,u)\mu(i,dq) + \Gamma\left(\int_{\mathcal{Q}(i)}\tilde{k}(i,q',u)\mu'(i,dq') - C\right) &\geq w_{i}, \ (j,u)\in\cS\times\Act, \nonumber\\
\mu, \ \mu', \ \nu, \ \nu' \in\PW.   & \ \nonumber
\end{align}

\bigskip

Here the Lagrange multiplier $\Gamma$ is unknown a priori. So one can use the following `primal-dual' scheme. Start with an initial guess for $\Gamma$, say $\Gamma_0 \geq 0$, and update it as follows. At step $n\geq 0$, solve the above linear programs for $\Gamma = \Gamma_n$. Let $\mu^n(\cdot, \cdot)$ be the optimal $\mu'(\cdot, \cdot)$ from the dual program and $y^n_\cdot(\cdot)$  the optimal $y_\cdot(\cdot)$ for the primal linear program, under $\Gamma = \Gamma_n$. Perform the iterate
$$\Gamma_{n+1} = \Gamma_n + a(n)\left(\int_{\widetilde{W}}\sum_u\tilde{k}(i,q',u)y_i^n(u)d\mu^n(i,dq') - C\right),$$
where $a(n) > 0$ is a stepsize sequence satisfying $a(n)\to 0, \ \sum_na(n) = \infty$.

\subsection{Reinforcement learning}
Reinforcement learning deals with data-driven algorithms for control. They are popular for situations when the system model is not known or is too messy, but adequate data, either real or simulated, online or offline, is available.  In most cases this is based on approximate dynamic programming. Thus these algorithms  usually  mimic classical iterative schemes for solution of dynamic programming equations. In fact they are usually stochastic approximation versions thereof. The development of reinforcement learning for ergodic risk-sensitive control, however, has been rather limited. We summarize it here. Much of the literature that talks of {\it risk-sensitive reinforcement learning} refers to other notions of risk arising from economics and finance. We do not consider them here, nor do we consider anything other than  the `ergodic' or time-asymptotic case that we have been considering so far. 
That is, we do not consider the finite horizon problem that has received some attention in literature \cite{Zaoran}.

The risk-sensitive Q-learning algorithm is inspired by the original Q-learning scheme for discounted cost \cite{Watkins}.
Consider DTCMC with finite state $\cS$ and finite action set $\Act$. Assume that the chain is irreducible and aperiodic under any stationary Markov policy. For notational ease, we replace $\lambda^{*,m}$ above by $\lambda^*$. Recall the risk-sensitive dynamic programming equation
\begin{equation}\label{DP}
V(i) = e^{-\lambda^*}\min_{u\in\Act}\left(e^{c(i,u)}\sum_{j\in\cS} P(j|i,u)V(j)\right), \ i\in\cS. 
\end{equation}
Setting $Q(i,u)\df$ the term in parenthesis on the right, we have a similar equation for $Q(\cdot,\cdot)$:
\begin{equation}\label{QDP}
Q(i,u) = \left(e^{c(i,u)-\lambda^*}\sum_{j\in\cS}P(j|i,u)\min_{u'}Q(j,u')\right), \ (i,u)\in\cS\times\Act, 
\end{equation}
where we have used the fact 
\begin{equation}\label{V-Q}
V(i) = e^{-\lambda^*}\min_{u\in\Act} Q(i,u) \ \forall \ i\in\cS.
\end{equation}
Replacing (\ref{DP}) by (\ref{QDP}) has increased the dimensionality from 
$|\cS|$ to $|\cS|\times|\Act|$, but the advantage is that the nonlinearity, that is, the `min' operator, is now \textit{inside} the conditional expectation. Since stochastic approximation  at its core is an averaging technique, this make the problem amenable to a model-agnostic, data-driven stochastic approximation algorithm.
Another advantage is that once you know $Q(\cdot,\cdot)$ exactly or approximately, dynamic programming tells us that the best control choice in state $i$ is $\Argmin(Q(i, \cdot))$. This does not require any knowledge of the model.

Note that both $V(\cdot)$ and $\lambda^*$ are unknowns. Taking cue from the average cost Q-learning \cite{ABB}, we have the algorithm inspired by {\it relative value iteration} for risk-sensitive control, based on a real or simulated run of a controlled Markov chain $(X_n,\zeta_n)$. This is given by:
\begin{align}
Q_{n+1}(i,u) = Q_n(i,u) + a(\nu(i,u,n))\Ind_{\{X_n = i, \zeta_n = u\}}\times\left(\frac{\min_{u'}Q(X_{n+1},u')}{Q_n(i_0,u_0)} - Q_n(i,u)\right), \ n \geq 0,\label{Qlearn}
\end{align}
where 
$$\nu(i,u,n) \df \sum_{m=0}^n \Ind_{\{X_m = i, \zeta_m = u\}}$$
is assumed to satisfy
$$\liminf_{n\to\infty}\frac{\nu(i,u,n)}{n} > 0 \ \mbox{a.s.} \ \forall \ i,u.
$$
That is, all state-control pairs are sampled  `comparably often', a.s. (This is a standard assumption for reinforcement learning algorithms.) The step-size sequence $\{a(n)\}$ satisfies
\begin{equation}
a(n) > 0, \ \sum_na(n) = \infty, \ \sum_na(n)^2 < \infty, \label{RM}
\end{equation}
plus some additional technical conditions in case of this specific (that is, fully asynchronous)  variant. Ignoring the technicalities due to asynchrony, the passage from what could have been a relative value iteration for $Q(\cdot,\cdot)$ and the foregoing is that one first replaces in the right hand side of (\ref{QDP})  the conditional expectation by an actual evaluation at a sample generated according to the conditional probability in question ($X_{n+1}$ in this case), replaces the unknown  $e^{\lambda^*}$ by $Q(i_0,u_0)$ for some fixed choice of $(i_0,u_0)$\footnote{Other choices for this `normalization' are possible.}, and then takes a convex combination of this with the previous iterate with weights $a(n), 1 - a(n)$, resp. Because we are considering a scheme based on a single run of the chain, we observe only a single transition at each time and therefore can update at time $n$ only the $(i,u)$-th component for which $X_n = i, \zeta_n = u$. Hence we multiply $a(n)$ by the indicator $\Ind_{\{X_n = i,\zeta_n=u\}}$ in the above, thus leaving the rest of the components unchanged.

Stochastic approximation theory then tells us that the iterates a.s.\ track the differential equation
$$
\dot{q}_t(i,u) = \frac{\sum_{j\in\cS} P(j|i,u)\min_{u'}q_t(j,u')}{q_t(i_0,u_0)} - q_t(i,u),
$$
which can be shown to converge to the solution $Q$ of (\ref{QDP}) for which the $(i_0,u_0)$th component is $e^{\lambda^*}$. Then so does $\{Q_n\}$, a.s. This can be viewed as a stochastic approximation version of the well known `power iteration' method of computational linear algebra, albeit for a nonlinear map. The analysis involves mapping the trajectories of this differential equation to those of a related differential equation given by
$$
\dot{q}_t'(i,u) = e^{-\lambda^*}\,\sum_{j\in\cS} P(j|i,u)\min_{u'}q_t'(j,u') - q_t'(i,u),
$$
which is easier to analyze. See \cite{MR1908528} for details.

One problem with the above `exact' Q-learning scheme, known as `tabular form' in machine learning literature, is that its dimensionality can be prohibitive. This prompts the use of a parametrized family of approximate $Q(\cdot,\cdot; \theta)$, where $\theta$ is a parameter vector of moderate dimensions, and then write a recursion for $\{\theta_n\}$ to learn the `best' $\theta$ in a suitable sense. One popular choice has been linear parametrization, i.e.,  a linear combination of a suitable choice of basis functions, both because of its ease and analytic tractability. One such scheme was studied in \cite{MR2464648}, albeit for policy evaluation, that is, for learning an approximate value function for a fixed randomized Markov policy, for which a rigorous theory is possible. (That linear function approximation may not work with the nonlinearity - the `min' operation - in place is a known fact even for simpler cost criteria.)  Interesting approximation error bounds for the eigenvalue $e^{\lambda^*}$ have been derived in \cite{Karm}. More recently, deep neural networks have been the favoured approximation architecture in other contexts, but they do not seem to have been explored in the risk-sensitive scenario.

Just as Q-learning is related to value iteration or relative value iteration as the case may be, another leading algorithm called the Actor-Critic algorithm is related to the policy iteration. Here we replace (\ref{Qlearn}) by
\begin{equation}
V_{n+1}(i) = V_n(i) + a(\nu(i,n))\Ind_{\{X_n = i\}}\left(\frac{V(X_{n+1})}{V_n(i_0)} - V_n(i,u)\right), \ n \geq 0,\label{Vlearn}
\end{equation}
with the randomized stationary Markov control policy
$\pi_n(i,u) \df \Prob(U_n=u|X_n= i)$ given recursively by a stochastic gradient scheme for the risk-sensitive cost. This recursion is performed with a different step-size sequence $\{b(n)\}$ which, in addition to satisfying the usual conditions (\ref{RM}), also satisfies $b(n) = \sorder(a(n))$, so that this iteration moves on a slower time scale. The net effect is that (\ref{Vlearn}) sees the latter as quasi-static, hence it can be analyzed by treating $\pi_n\approx$ constant,
leading to the conclusion that it is `essentially' a policy evaluation scheme that tracks the value function for the \textit{constant} policy $\pi \approx \pi_n$. That is, denoting by $V_\pi$ the value function for a fixed stationary policy $\pi$, we have $V_n - V_{\pi_n} \to 0$ a.s. This emulates the policy evaluation component of policy iteration. The $\{\pi_n\}$-iterate, performing gradient descent as though $V_n$ is a legitimate surrogate for $V_{\pi_n}$ (which it is, as argued above) emulates the optimization step of policy iteration. We omit the details of the latter, suffice to say that it is based on a sensitivity formula for risk-sensitive cost with respect to a parameter. See \cite{MR2021952} for details. One limitation of this work is that the optimization component works only if you update $\pi_n$ directly, not its parametrized approximation, because such variants require model knowledge for their implementation.  Recently a policy gradient scheme based on updates only at successive visits to a privileged state has been proposed as a workaround \cite{MMRS}.
For (\ref{Vlearn}), however, its `essentially linear' nature allows for justifiable use of linear parametrization, that is, as in \cite{MR2464648}. It was recently observed in \cite{BBG}, albeit for a different cost, that interchange of the fast and slow time scales in Actor-Critic algorithm leads to a new algorithm that emulates value iteration, dubbed `Critic-Actor algorithm' in {\it ibid}. This also applies to risk-sensitive problem. One loses, however, the legitimacy of linear function approximation.

Given the thin list of references here, it is clear that this remains a wide open area for further research.

\subsection*{Acknowledgments}
This research of Anup Biswas was supported in part by a SwarnaJayanti fellowship  DST/SJF/MSA-01/2019-20. The research of Vivek Borkar was supported in part by a S.\ S.\ Bhatnagar Fellowship.

\bibliography{RISK-main}

\end{document}